\documentclass[compress,final,3p,times,12pt]{elsarticle}
\usepackage{amsfonts}
\usepackage{graphics}
\usepackage{amssymb}
\usepackage{amsmath}
\usepackage{color, epstopdf}
\usepackage{mathrsfs}

\usepackage{graphicx, color, epstopdf}
\usepackage{threeparttable}
\usepackage{float}
\usepackage{array}
\usepackage{subcaption}

\numberwithin{equation}{section}
 \numberwithin{Lem}{section}
 \numberwithin{Defi}{section}
 \numberwithin{Theo}{section}
 \numberwithin{Rem}{section}
  \numberwithin{Coro}{section}
  \numberwithin{Fig}{section}

\journal{}

\begin{document}

\begin{frontmatter}



\title{Most probable dynamics of a genetic regulatory network under stable L\'{e}vy noise} \tnotetext[label1]{This work is
 supported by NSFC (Grant Nos. 11531006, 11371367, 11271290, 11771062) and the Fundamental Research Funds for the Central Universities', HUST: 0118011076.}
\author{Xiaoli Chen\fnref{addr1}}
\ead{xlchen@hust.edu.cn}
\author{Fengyan Wu\fnref{addr3}}
\ead{fywuhust@163.com}\cortext[cor1]{Corresponding author}
\author{Jinqiao Duan\fnref{addr5}}
\ead{duan@iit.edu}
\author{J\"{u}rgen Kurths\fnref{addr6,addr7}}
\ead{kurths@pik-potsdam.de}
\author{Xiaofan Li\fnref{addr5}}
\ead{lix@iit.edu}

\address[addr1]{Center for Mathematical Sciences \& School of Mathematics and Statistics, Huazhong
University of Science and Technology, Wuhan 430074,  China}
\address[addr3]{College of Mathematics and Statistics, Chongqing University, Chongqing 401331, China}%
\address[addr5]{Department of Applied Mathematics, Illinois Institute of Technology, Chicago, IL 60616, USA}
\address[addr6]{ Research Domain on Transdisciplinary Concepts and Methods,
Potsdam Institute for Climate Impact Research, PO Box 60 12 03, 14412 Potsdam, Germany}
\address[addr7]{Department of Physics, Humboldt University of Berlin, Newtonstrate 15, 12489 Berlin, Germany}
\begin{abstract}
{Numerous studies have demonstrated the important role of noise in the dynamical behaviour of a complex system. The most probable trajectories of nonlinear systems under the influence of Gaussian noise have recently been studied already. However, there has been only a few works that examine how most probable trajectories in the two-dimensional system (MeKS network) are influenced under non-Gaussian stable L\'{e}vy noise. Therefore, we discuss the most probable trajectories of a two-dimensional model depicting the competence behaviour in \emph{B. subtilis} under the influence of stable L\'{e}vy noise. On the basis of the Fokker-Planck equation, we describe the noise-induced most probable trajectories of the MeKS network from the low ComK protein concentration (vegetative state) to the high ComK protein concentration (competence state) under stable L\'{e}vy noise. We demonstrate choices of the non-Gaussianity index $\alpha$ and the noise intensity $\epsilon$ which generate the ComK protein escape from the low concentration to the high concentration. We also reveal the optimal combination of both parameters $\alpha$ and $\epsilon$ making the tipping time shortest. Moreover, we find that different initial concentrations around the low ComK protein concentration evolve to a metastable state, and provide the optimal $\alpha$ and $\epsilon$ such that the distance between the deterministic competence state and the metastable state is smallest.}
\end{abstract}


\end{frontmatter}

\section{\label{sec:level1}Introduction}
\subsection{\label{sec:A} Background}
A number of studies in different fields were concerned about the role of noise playing in dynamical systems, since noise may induce various delicate effects \cite{oksendal,Applebaum,DuanBook,Bressloff,Kaern05,Raj08}. Noise-induced transition phenomena were regarded in various dynamical systems including gene regulatory networks \cite{Suel06,Paulcf4}, oscillators and electronic transport in physics \cite{Paulcf2,Paulcf3,Paulcf5, Swain1,Franovi}, chemical reacting systems \cite{Paulcf1,larrisacf16}, cancer cell proliferation and epidemiological models in pathology \cite{Paulcf7,Paulcf8}.

These applications have been considered for Guassian noise. But when fluctuations are present in certain events, such as burst-like events, Gaussian noise is not appropriate. In this case, it is more appropriate to model the random
fluctuations by a non-Gaussian L\'{e}vy motion with heavy tails and bursting sample paths. Non-Gaussion L\'{e}vy noise has found a lot of applications in different areas. Perc applied the non-Gaussian L\'{e}vy noise in economy and studied the impact of an important class of external perturbations on the evolution of cooperation in a spatially extended prisoner's dilemma game \cite{Perc1,Perc2}.
L\'{e}vy noise was also associated with saving energy in animals, e.g. albatross flight patterns, as reviewed in \cite{Trenchard}. Zheng et. al considered how L\'{e}vy noise may affect gene regulatory networks in \cite{yayun16}. For further applications about non-Gaussion L\'{e}vy noise cf. \cite{hui1,Cai,la16,hui2}.

Gene expression is a noisy process \cite{Swain1,Marinov}. There are various studies about Gaussian noise in gene regulation \cite{Simpson,Wells,Gui,huicf9,huicf10}. However, the transcription of a gene can be a discontinuous or burst-like event,
and mRNA is synthesized in intermittent but intense pulses or bursts \cite{yayun16,hui1,burst1,burst2,burst3}. Under this circumstance, Gaussian noise is not proper to model this phenomenon, but a stochastic process with discontinuous and heavy tail distribution, i.e., non-Gaussian L\'{e}vy noise \cite{huicfxuy,Jia,Kumar,Holloway,levyEvid,smolen98} appears more appropriate in describing the fluctuations in gene regulation.

In this study, we devote to studying the transitions in a two-dimensional genetic regulatory model driven by non-Gaussian $\alpha$-stable L\'{e}vy noises, which describe fluctuations with features such as heavy tails and jumps. The corresponding genetic regulatory model was established by S\"{u}el et al. \cite{Suel06}. Comparing with the experimental data, they demonstrated  the reasonability and validity of the genetic regulatory model and  considered a set of excursion trajectories in an excitable case under Gaussian noise. Instead of using the mean first exit time and first escape probability to study the bistable case of this model driven by $\alpha$-stable L\'{e}vy noise \cite{wu17}, we examine the most probable trajectories rather than stochastic trajectory sample paths to characterize the dynamical behaviors of the model. The nonlocal Fokker-Planck equation of stochastic systems with L\'evy noise can be numerically simulated \cite{Gaoting16,Gaoting14,Li1,Li2,Li3,Chen1,Chen2}. The noise induced most probable trajectories, which describe the MeKS network from the low ComK protein concentration (vegetative state) to the high ComK protein concentration (competence state), are obtained by examining the solution of the nonlocal Fokker-Planck equation. We use here numerical methods to solve the nonlocal Fokker-Plank equation, (for more details see Appendixes). We present choices of the non-Gaussianity index $\alpha$ and the noise intensity $\epsilon$ which induce the ComK protein transfer from the low concentration to the high concentration. Then we will analyse how different initial concentrations around the low ComK protein concentration affect the metastable state. Furthermore, we determine the optimal $\alpha$ and $\epsilon$ such that the distance between the deterministic competence state and the metastable state is smallest.

\subsection{\label{sec:B} A stochastic model for the MeKS network}
To investigate the dynamics of competence induction, S\"{u}el et al. \cite{Suel06} constructed a MeKS network:
\begin{equation}  \label{model}
\begin{split}
\tfrac{dk}{dt}&={{a}_{k}}+\tfrac{{{b}_{k}}{{k}^{n}}}{k_{0}^{n}+{{k}^{n}}}-\tfrac{k}{1+k+s},\\
 \tfrac{ds}{dt}&=\tfrac{{{b}_{s}}}{1+{{(k/{{k}_{1}})}^{p}}}-\tfrac{s}{1+k+s}.
 \end{split}
\end{equation}
The symbols $k$ and $s$ represent the concentration of the ComK and ComS proteins, respectively. The parameters $a_k$ and $b_k$ stand for the basal and fully activated rates of ComK protein production, respectively. The parameter $k_0$ denote the concentration of the ComK protein needed for a $50\%$ activation. The Hill coefficients $n$ and $p$ are the cooperativities of the ComK protein auto-activation and ComS protein repression, respectively. The maximal expression rate of the ComS protein is $b_s$, and when $k=k_1$ the ComS protein attains its half-maximal rate. We choose suitable parameters \cite{Suel06}: $a_k\!=\!0.004,~b_k\!=\!0.14,~b_s\!=\!0.68,~k_0\!=\!0.2,~k_l\!=\!0.222,~n\!=\!2,~p\!=\!5$. In this case, the MeKS network has three equilibria: two stable equilibria  (0.015262, 2.1574) (nodal sink corresponding to the low vegetative  state), (0.15732,1.5781) (spiral sink corresponding to the high competence state) and one unstable equilibrium (0.08568,2.2469)(saddle point).

\begin{figure}[H]
\hspace{4cm}
        \includegraphics[width=0.5\textwidth]{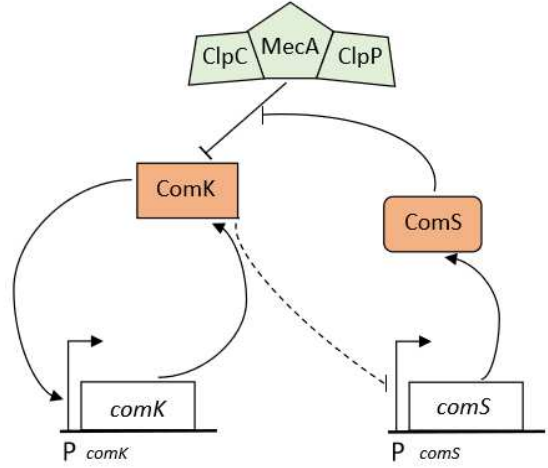}
    \caption{\textbf{The core competence circuit in \emph{B. subtilis} (MeKS).} $P_{comK} (P_{comS})$ denotes the promoter of $\emph{comK} ~ (\emph{comS})$ gene. There is an autoregulatory positive feedback loop regulated by the ComK protein, and a negative feedback loop, in which the ComK protein indirectly represses the expression of the comS gene and the ComS protein competitively interferes with the proteolytic degradation of ComK protein.}
  \label{f1}
\end{figure}

The core competence genetic circuit of the MeKS network is shown in FIG. \ref{f1}. In \emph{B. subtilis} population, mostly, the ComK protein is expressed by the \emph{comK} gene at a relative low level. Thus, the ComK protein concentration is low, which corresponds to the vegetative state. Under circumstance of nutrient limitation, a small part of \emph{B. subtilis} differentiate into the competence state, in which the ComK protein concentration is high. In this case, the cell can take in extracellular DNA from the environment \cite{Suel06,Dubnau99,Grossman95,Suel09,Suel15}. It was found that the key transcription factor, the ComK protein, can activate the transcription of several genes which is necessary for the competence development. When the ComK protein concentration is high, then the competence develops \cite{Pal13cf20,Pal13cf12}, from which we deduce that the ComK protein occupies a core position in the competence-signal-transduction network.

As discussed in the preceding subsection, gene expression in cells is subject to random fluctuations \cite{Swain1,Marinov}, and such fluctuations are under real conditions typically non-Gaussian noise. Hence we consider the MeKS network under the influence of $\alpha$-stable L\'{e}vy motion, which can be written as the following stochastic differential equations:
\begin{equation}  \label{stomodel}
 \begin{split}
  dk&=f_1(k,s)dt+\epsilon_k dL_t^1,  \\
  ds&=f_2(k,s)dt+\epsilon_s dL_t^2,
 \end{split}
\end{equation}
where $f_1(k,s)={{a}_{k}}+\tfrac{{{b}_{k}}{{k}^{n}}}{k_{0}^{n}+{{k}^{n}}}-\tfrac{k}{1+k+s},
 f_2(k,s)=\tfrac{{{b}_{s}}}{1+{{(k/{{k}_{1}})}^{p}}}-\tfrac{s}{1+k+s}.$
Here $(\epsilon_k,\epsilon_s)$ denote the diffusion coefficients,
 and $(L_t^1, L_t^2)$ is the two-dimensional L\'evy motion with independent scalar symmetric $\alpha$-stable L\'evy motions which have the same generating triplet $(0, 0, \nu_\alpha)$.  Here the  jump measure is a Borel measure in $\mathbb{R}\!\setminus\!\{0\}$ and defined by
  $$\nu_\alpha({\rm d}k)=C_\alpha|k|^{-(1+\alpha)}\, {\rm d}k,$$
where $\displaystyle{C_{\alpha} =
\frac{\alpha}{2^{1-\alpha}\sqrt{\pi}}
\frac{\Gamma(\frac{1+\alpha}{2})}{\Gamma(1-\frac{\alpha}{2})}}$, $\alpha \in (0, 2)$ is the
non-Gaussianity index. For more information about $\alpha$-stable L\'evy motions,
see \cite{Applebaum,DuanBook}.

 In the following discussions, we use the scale transformation $k'=10k,~s'=2s$. For the noise intensity, we set $\epsilon=\epsilon_k=\epsilon_s$.
 After this transformation, the two stable equilibria become $(k_{-},s_{-})=(0.15262,4.3148)$ (the low concentration state) and
$(k_{+},s_{+})=(1.5732,  3.1562)$ (the high concentration state), and the unstable equilibrium becomes $(k_{u},s_{u})=(0.8568,4.4938)$. We define $D_1\!=\!(0,k_u)$ as the low concentration region. If the most probable trajectories reach the high concentration region, where the ComK protein concentration is high, the competence develops.

This paper is organized as follows. In Section \ref{Method}, we present the method how to compute the most probable trajectories. In Section \ref{Result}, we discuss the most probable dynamics of the MeKS network with geometric tools (most probable trajectories) under stable L\'{e}vy noise. We finish this paper with conclusions and discussions in Section \ref{Conclusion}.
\section{\textbf{Method}}\label{Method}
\begin{figure}
    \begin{subfigure}[b]{0.4\textwidth}
    \leftline{~~~~~~~\tiny\textbf{(a)}}
        \includegraphics[width=\textwidth]{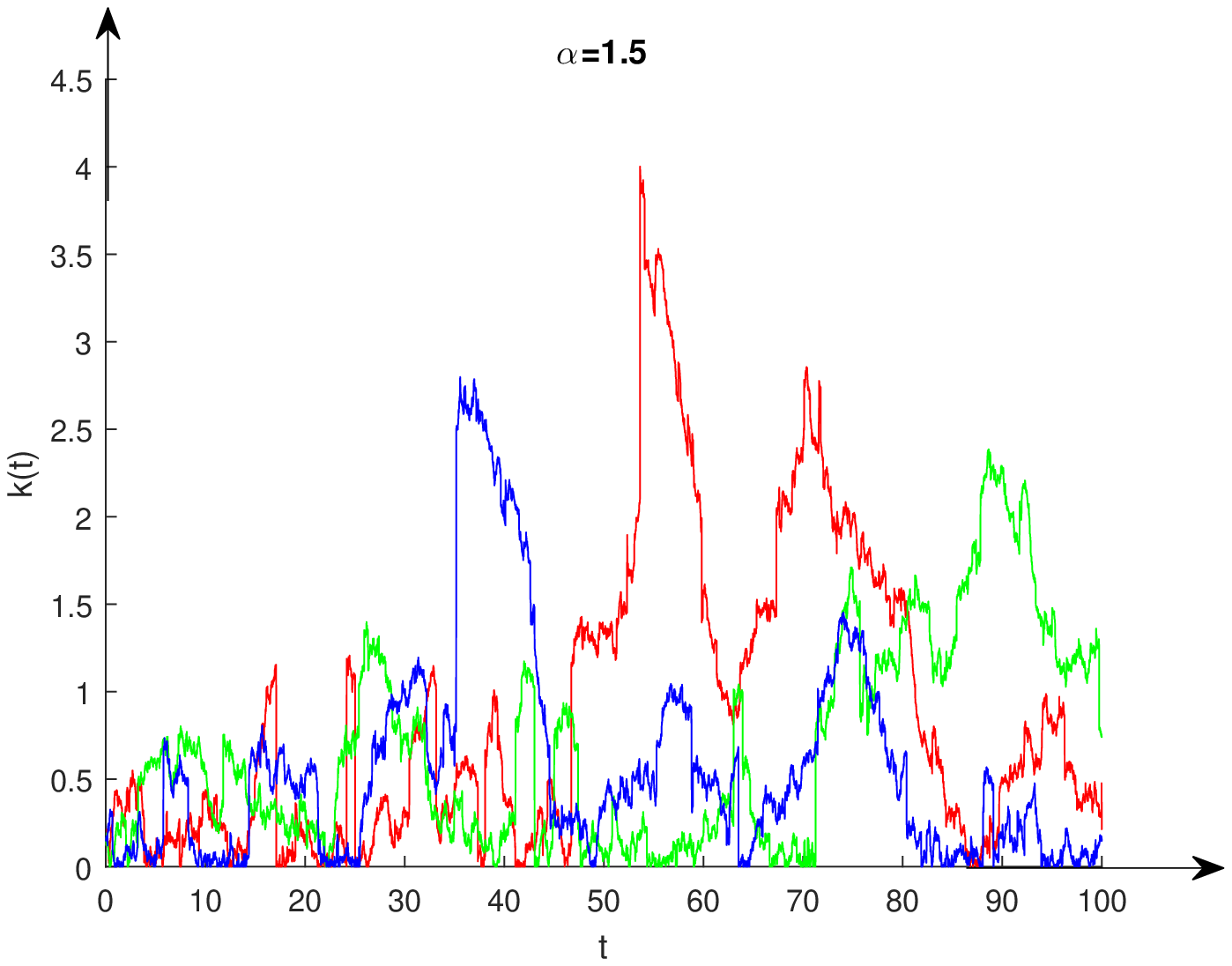}
    \end{subfigure}
   \begin{subfigure}[b]{0.4\textwidth}
   \leftline{~~~~~~~\tiny\textbf{(b)}}
        \includegraphics[width=\textwidth]{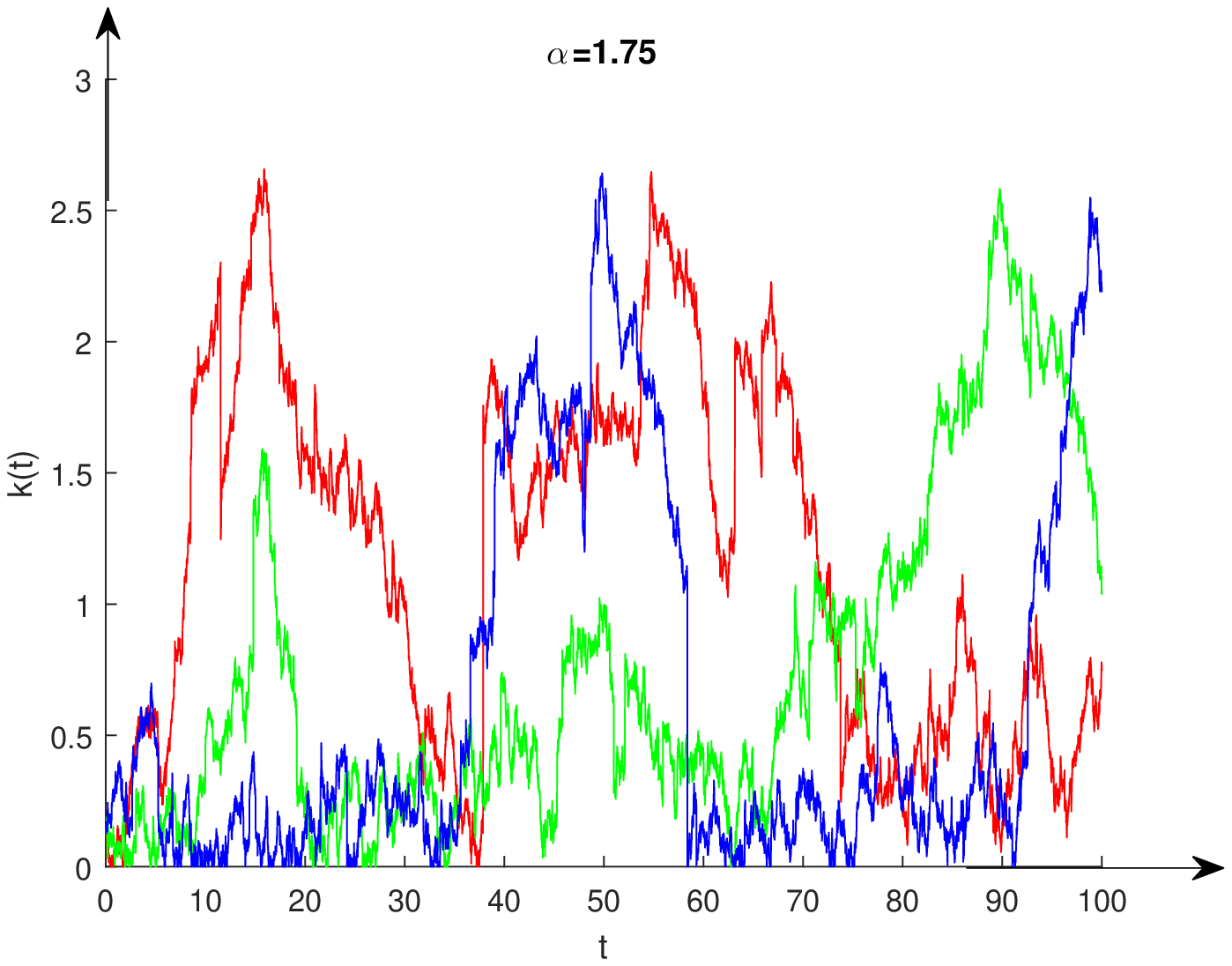}
    \end{subfigure}
\caption{\textbf{The stochastic trajectories of equation \eqref{stomodel} for ComK protein at $T=100$ for the noise intensity $\epsilon=0.15$. }
(a) The stochastic trajectories for $\alpha=1.5$.
(b) The stochastic trajectories for  $\alpha=1.75$.}
\label{f2}
\end{figure}

When it comes to trajectories for a stochastic dynamical system, there is an apparent option of consideration: the almost sure trajectories \cite{DuanBook,zhuan16}. However, as we know, the almost sure trajectories, i.e. the sample trajectories (FIG. \ref{f2}), which look like ``noodles" in the phase plane, could hardly provide helpful information for understanding the system's dynamics. Furthermore, to acquire comprehension of a stochastic dynamical phenomenon, stationary probability density functions for solution trajectories have been widely utilised \cite{zhuancf11,zhuancf9}.

Notice that each sample trajectory is an ``outcome'' of a trajectory of a stochastic dynamical system, from a given initial state. Now, the question arises, which trajectory does the system most probably to evolve in phase plane? Our geometric tool, the most probable trajectories \cite{DuanBook,zhuan16}, will tackle this problem. As is known, the solution of a stochastic system \eqref{stomodel} is a stochastic process, and the probability density function of the stochastic process is governed by the corresponding Fokker-Planck equation \cite{Gaoting16,zhuancf9}. On the basis of it, we can get the most probable trajectory by computing the maximizer of the probability density function $p(k,s,t)$ at every time $t$. This offers geometric representations of most probable trajectories of a stochastic dynamic system \cite{zhuan16}.

\begin{figure}[htp]
\begin{minipage}[]{0.5 \textwidth}
\centerline{\includegraphics[width=6cm,height=4cm]{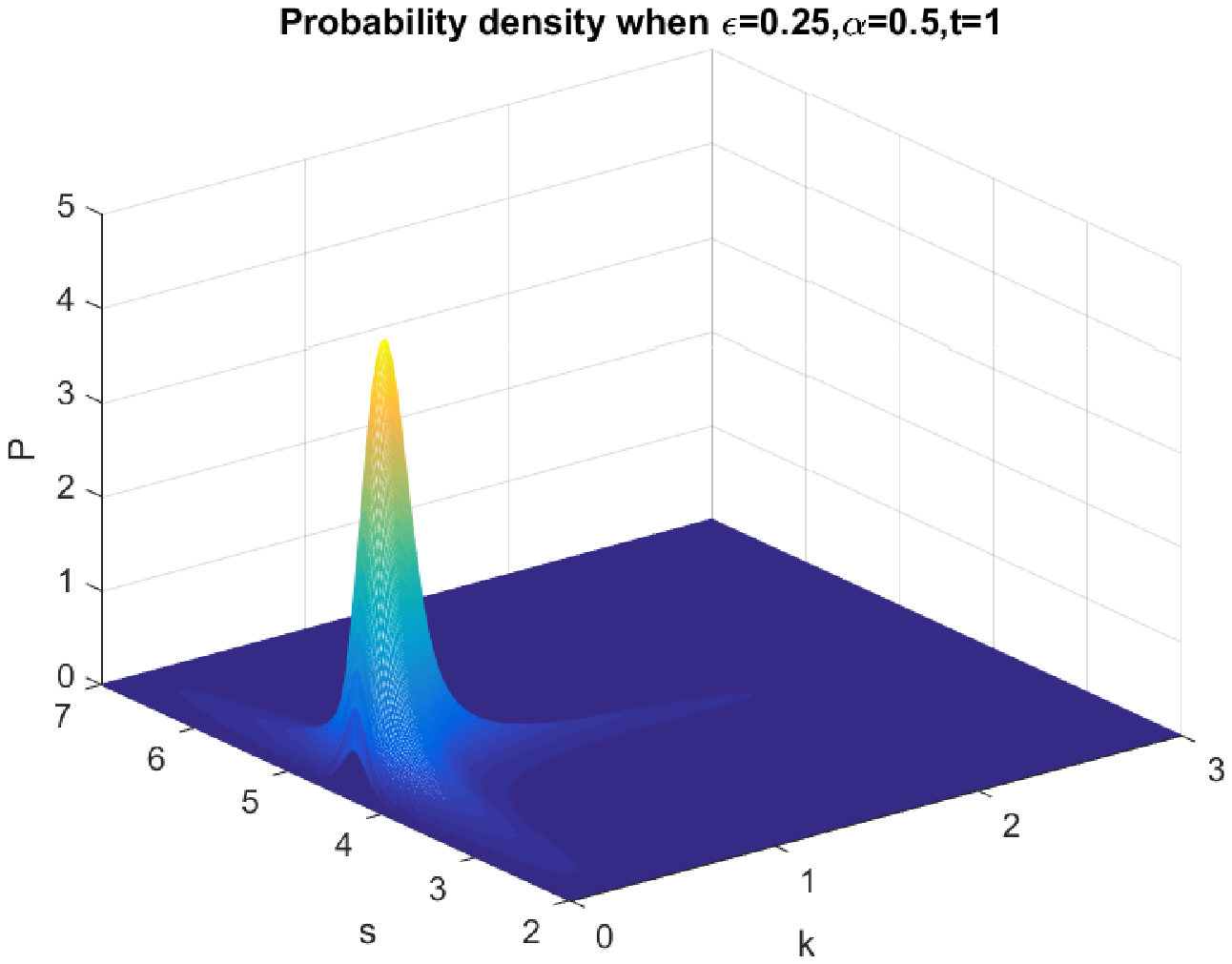}}
\centerline{(a) $t=1$}
\end{minipage}
\hfill
\begin{minipage}[]{0.5 \textwidth}
\centerline{\includegraphics[width=6cm,height=4cm]{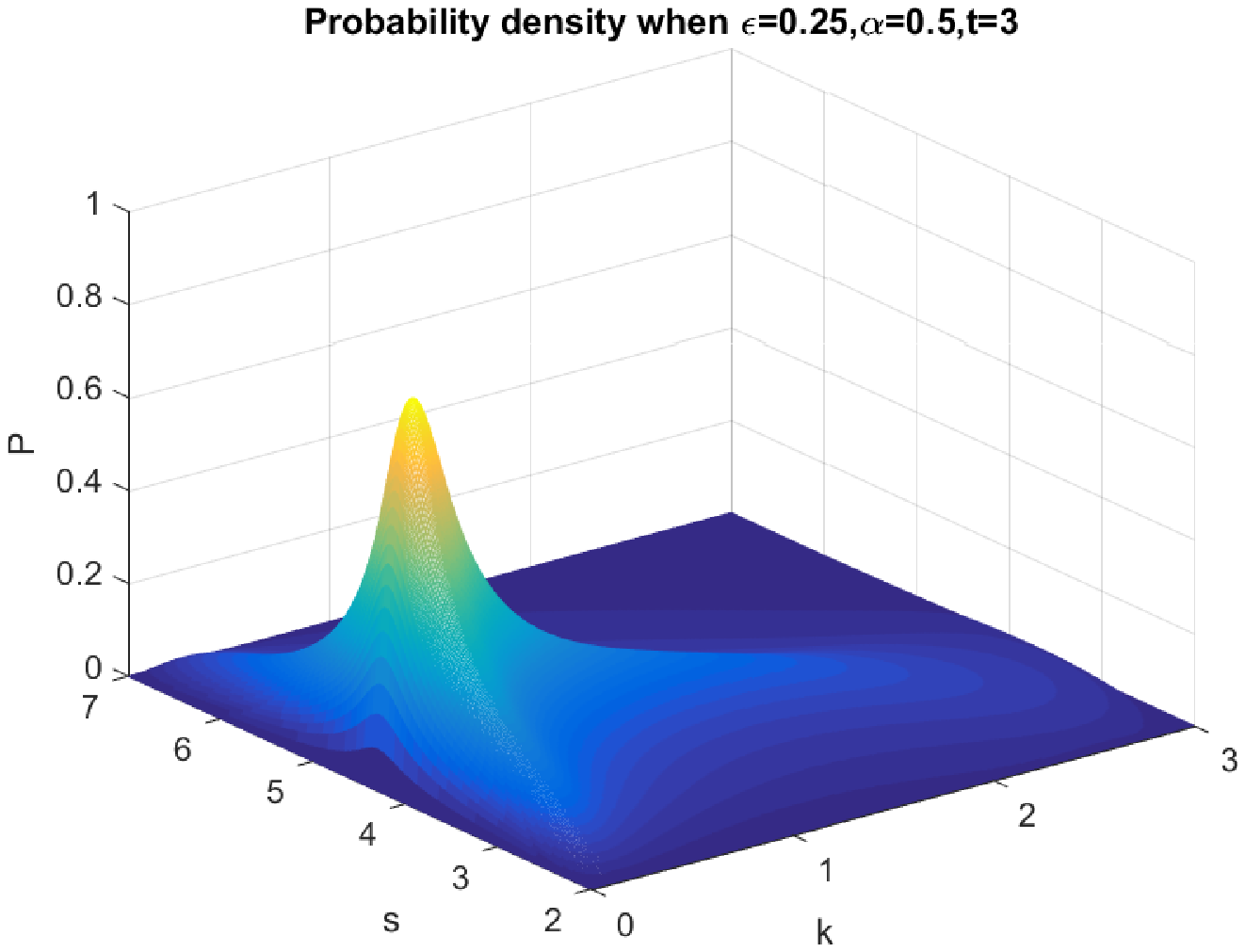}}
\centerline{(b) $t=3$}
\end{minipage}

\begin{minipage}[]{0.5 \textwidth}
\centerline{\includegraphics[width=6cm,height=4cm]{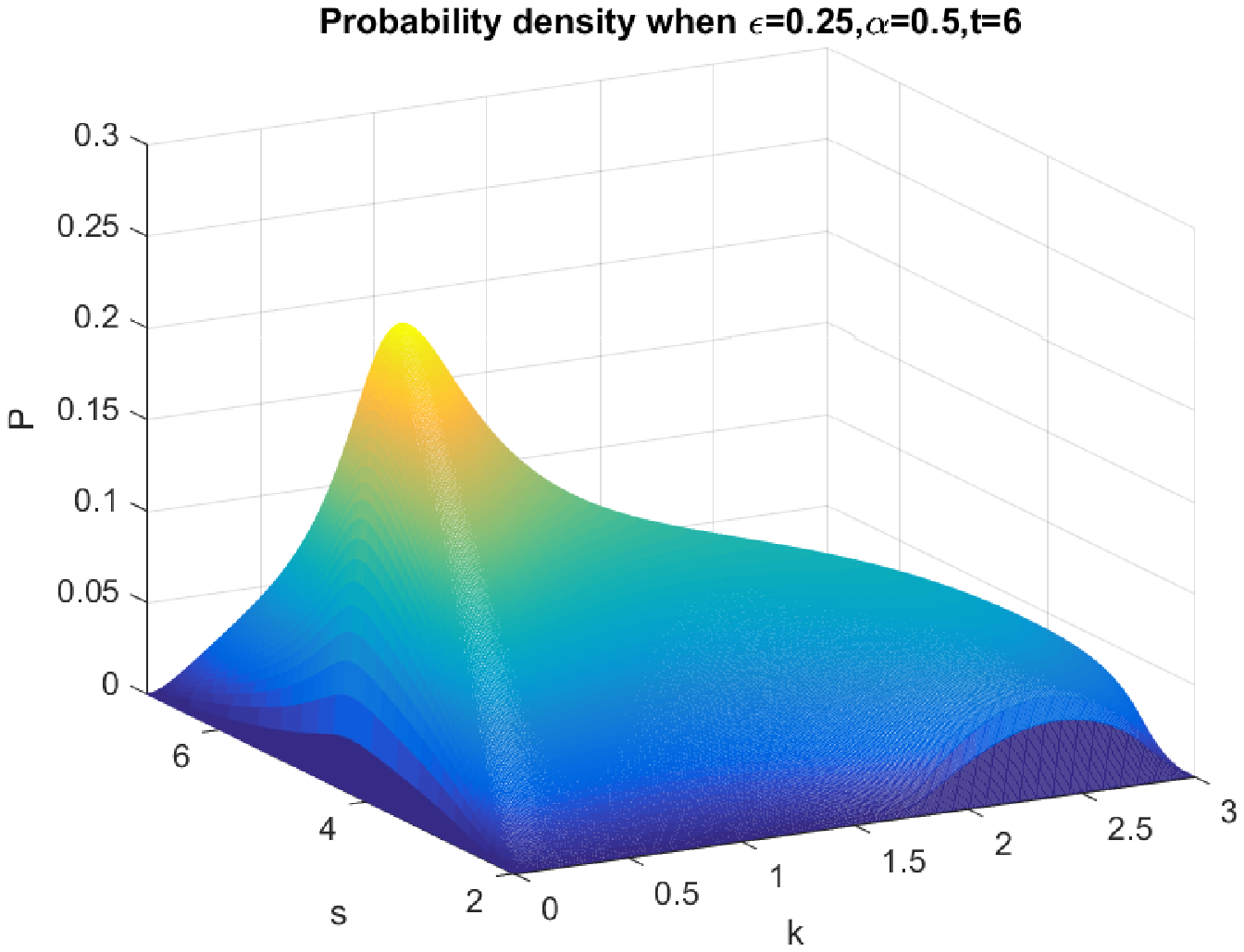}}
\centerline{(c) $t=6$}
\end{minipage}
\hfill
\begin{minipage}[]{0.5 \textwidth}
\centerline{\includegraphics[width=6cm,height=4cm]{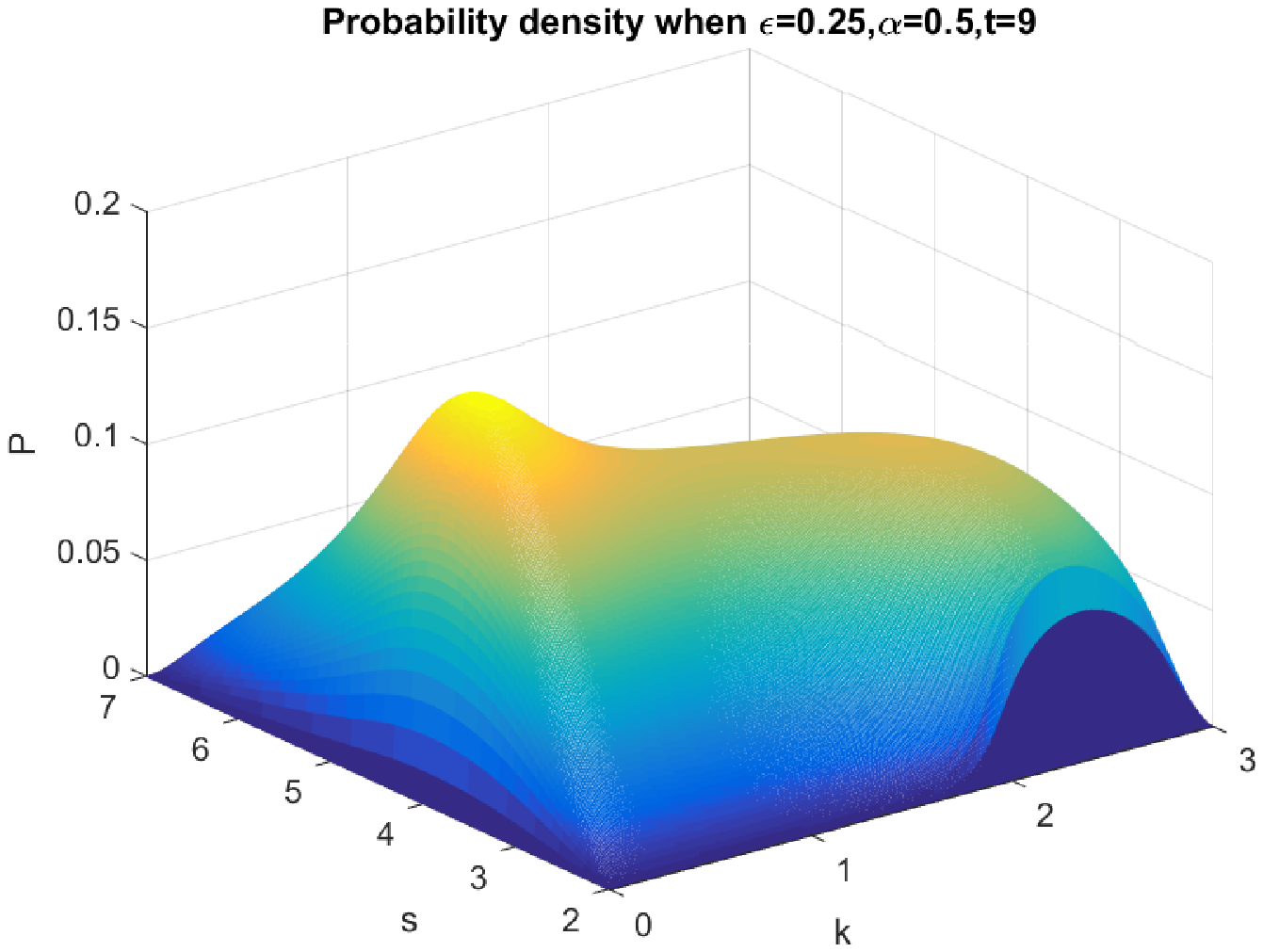}}
\centerline{(d) $t=9$}
\end{minipage}

\begin{minipage}[]{0.5 \textwidth}
\centerline{\includegraphics[width=6cm,height=4cm]{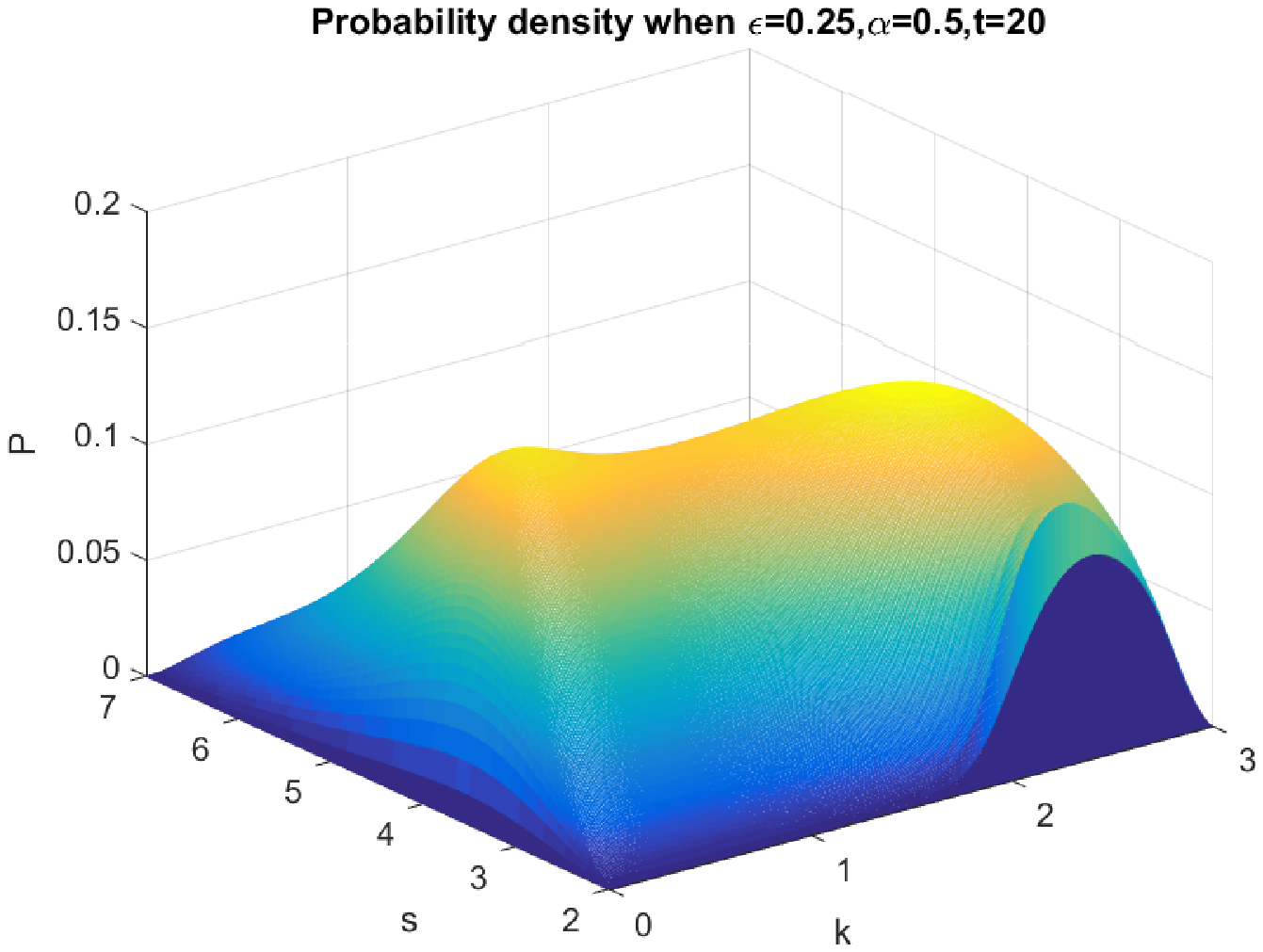}}
\centerline{(e) $t=20$}
\end{minipage}
\hfill
\begin{minipage}[]{0.5 \textwidth}
\centerline{\includegraphics[width=6cm,height=4cm]{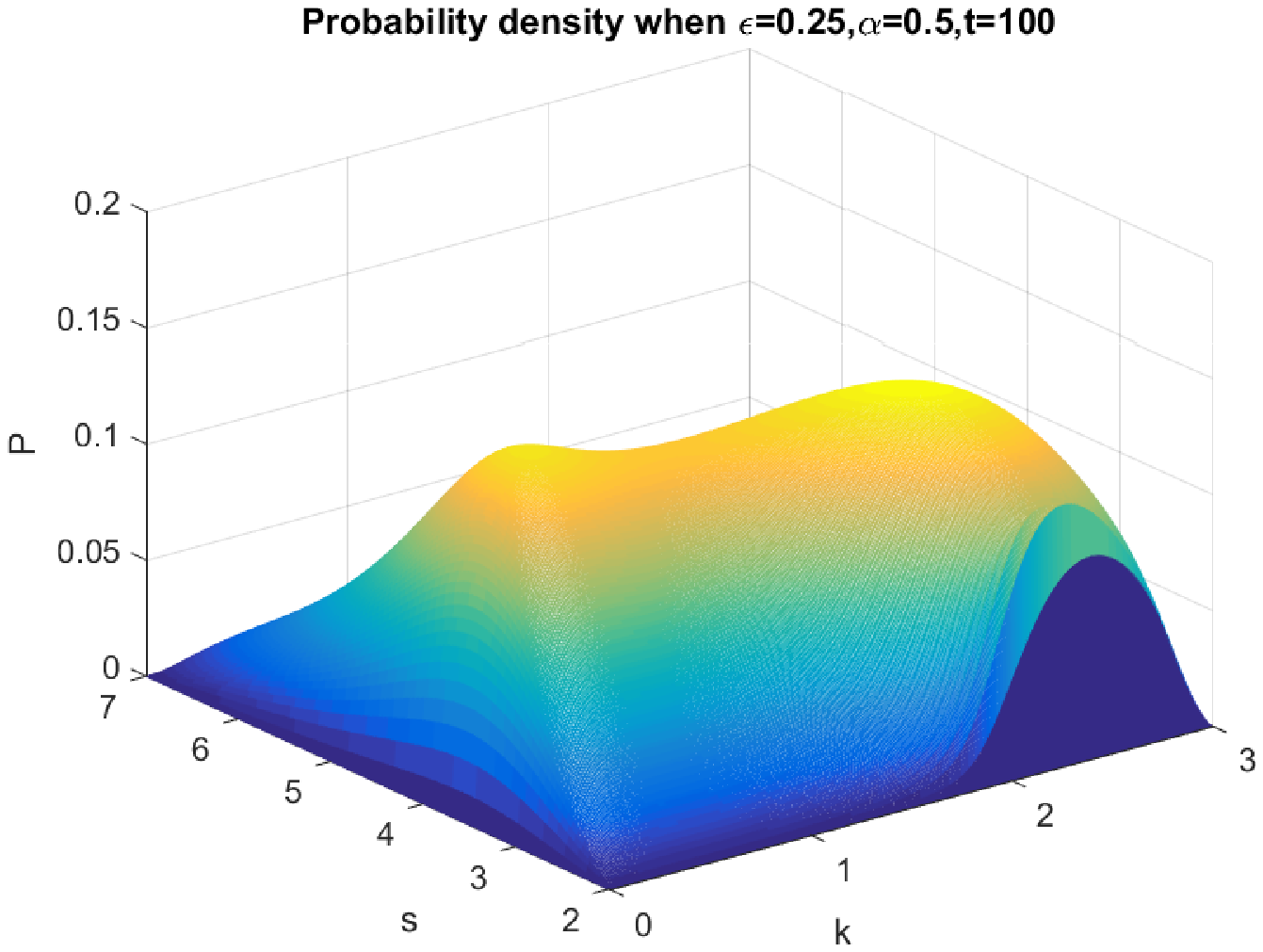}}
\centerline{(f) $t=100$}
\end{minipage}
\vfill
\caption{\textbf{Plots of the probability density function at various time: for $\epsilon=0.25$, $\alpha=0.5$ and the initial state $(x_0,y_0)=(0.15262,4.3148)$} .}
\label{f3}
\end{figure}

First we derive the generator and the Fokker-Planck equation of system \eqref{stomodel} (for more details see Section \ref{A}). Then we use the finite difference method to solve the Fokker-Planck equation (more details are given in Section \ref{B}). Plots of the so calculated probability density functions are shown in FIG. \ref{f3}. We set the low concentration state $(k_{-},s_{-})$ as the initial point of probability density function. FIG. \ref{f3} shows that the peaks transfer from one stable point to another with increasing time. When time is small $(t=1)$, there is one peak around the low concentration. As time goes on, there are two peaks for the probability density, one peak is transferred to the competence region, where the ComK protein concentration is high, the other peak still stays around the vegetative region. For $t=20$, the probability density function has attained its stationary solution, finally it stays at the competence state (the high ComK concentration). For a fixed time $t_i$, we compute the maximum of $p(k,s,t_i)$ to get the position $(k_i,s_i)$. We connect this series of $\{(k_i,s_i),~i=1,~2,\cdot\cdot\cdot\}$, to get the most probable trajectory (from a given initial state). The time instants $t_i$'s need to be taken close enough, in order to get a reasonable approximation of the most probable trajectory.

Through computing the most probable trajectories, we analyze how the ComK protein transfers from the low concentration regime to the high one, with different noise parameters. Then we want to know when the most probable trajectories go
over the saddle point. We define this time as the \emph{tipping time}. Moreover, we will
determine the choices of $\alpha$ and $\epsilon$ to benefit the development to competence.

\section{Results}\label{Result}
In this section, we analyse how the noise intensity $\epsilon$ and non-Gaussianity index $\alpha$   affect   the most probable trajectories. As seen from FIG. \ref{f3}, the probability density function attains its stationary state as time goes on. In the following, we choose $T=100$ as the computational terminal time.
\begin{figure}[htp]
\begin{minipage}[]{0.5 \textwidth}
 \leftline{~~~~~~~\tiny\textbf{(a1)}}
\centerline{\includegraphics[width=6cm,height=4cm]{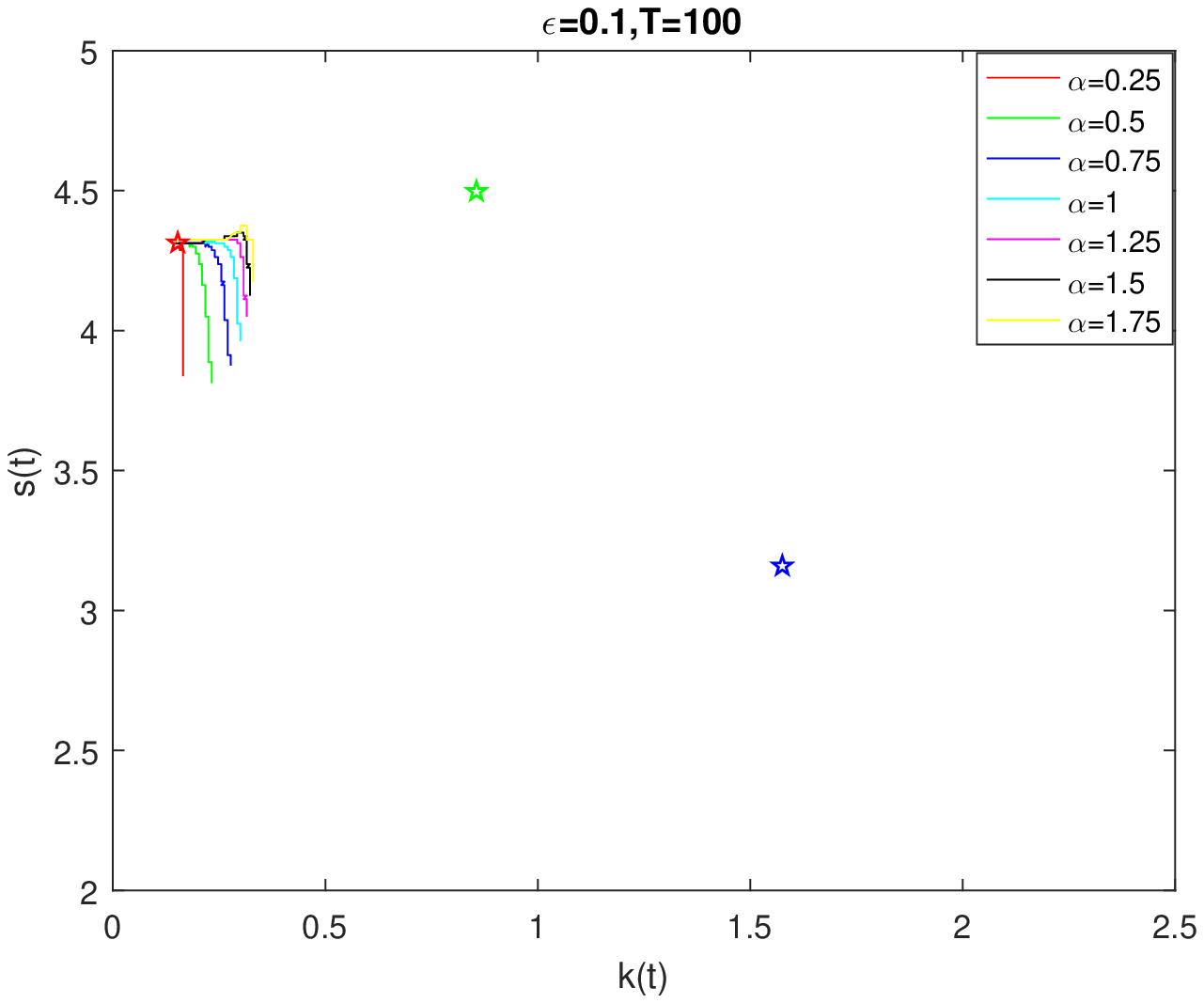}}
\end{minipage}
\hfill
\begin{minipage}[]{0.5 \textwidth}
 \leftline{~~~~~~~\tiny\textbf{(a2)}}
\centerline{\includegraphics[width=6cm,height=4cm]{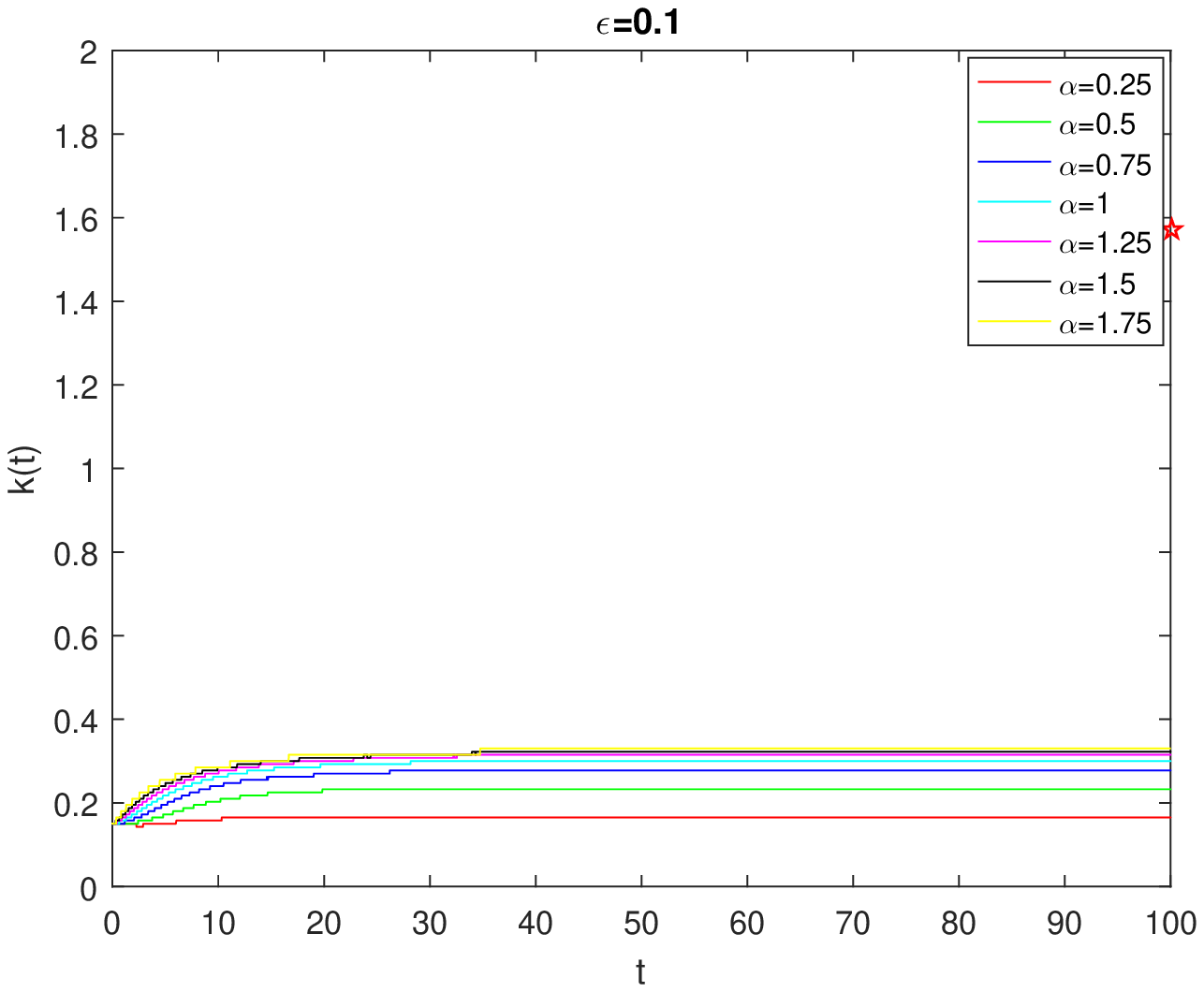}}
\end{minipage}

\begin{minipage}[]{0.5 \textwidth}
 \leftline{~~~~~~~\tiny\textbf{(b1)}}
\centerline{\includegraphics[width=6cm,height=4cm]{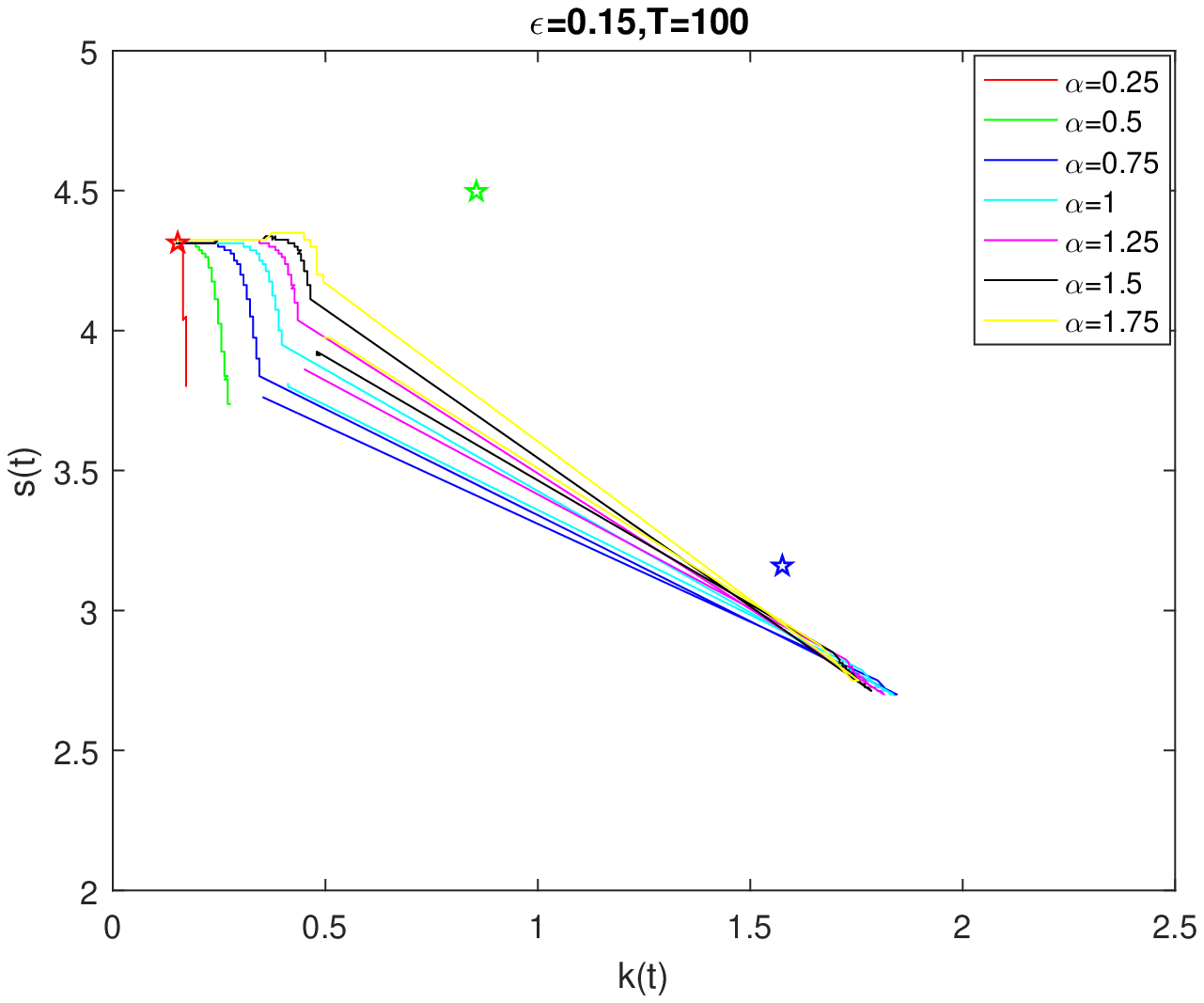}}
\end{minipage}
\hfill
\begin{minipage}[]{0.5 \textwidth}
 \leftline{~~~~~~~\tiny\textbf{(b2)}}
\centerline{\includegraphics[width=6cm,height=4cm]{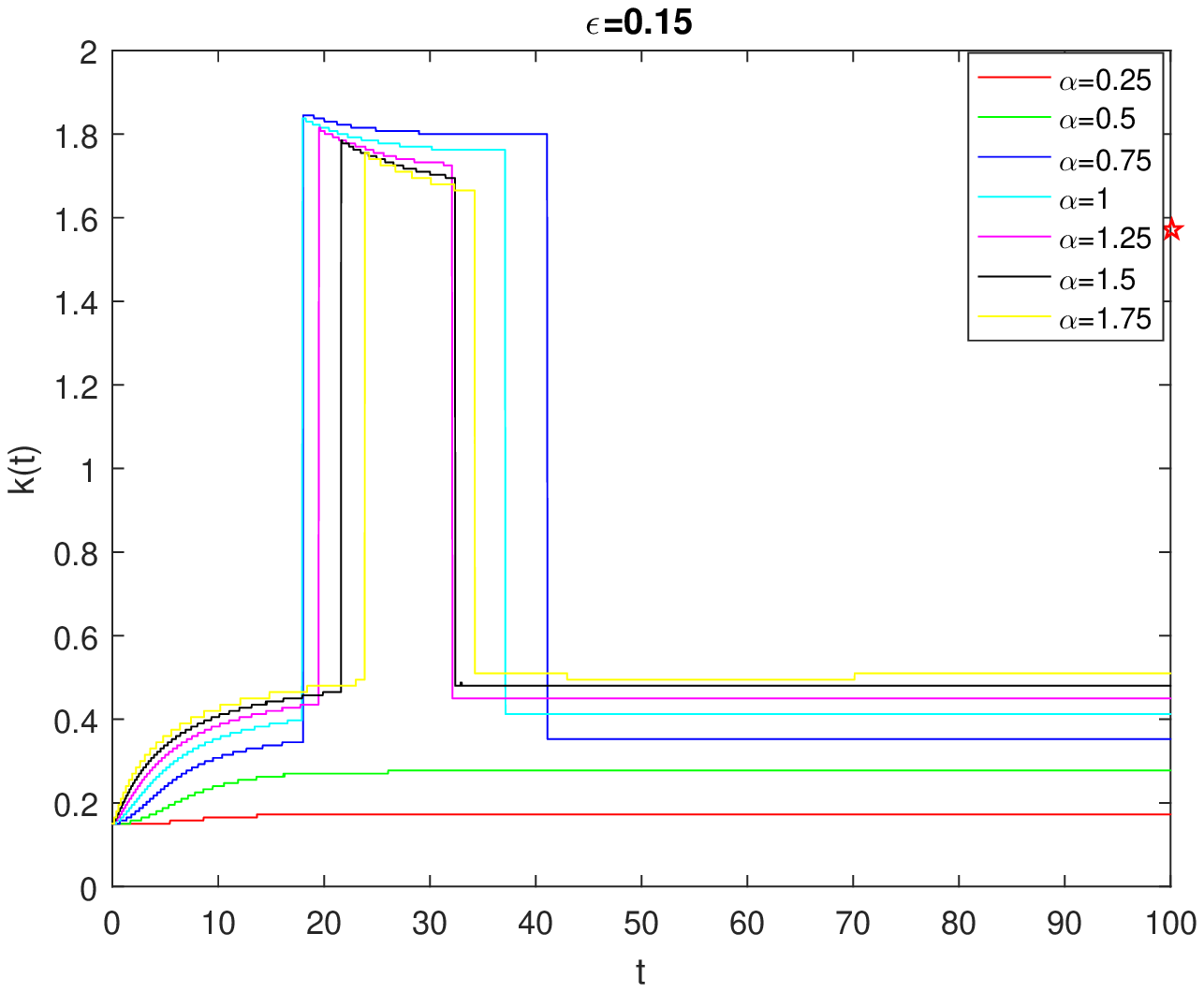}}
\end{minipage}

\begin{minipage}[]{0.5 \textwidth}
 \leftline{~~~~~~~\tiny\textbf{(c1)}}
\centerline{\includegraphics[width=6cm,height=4cm]{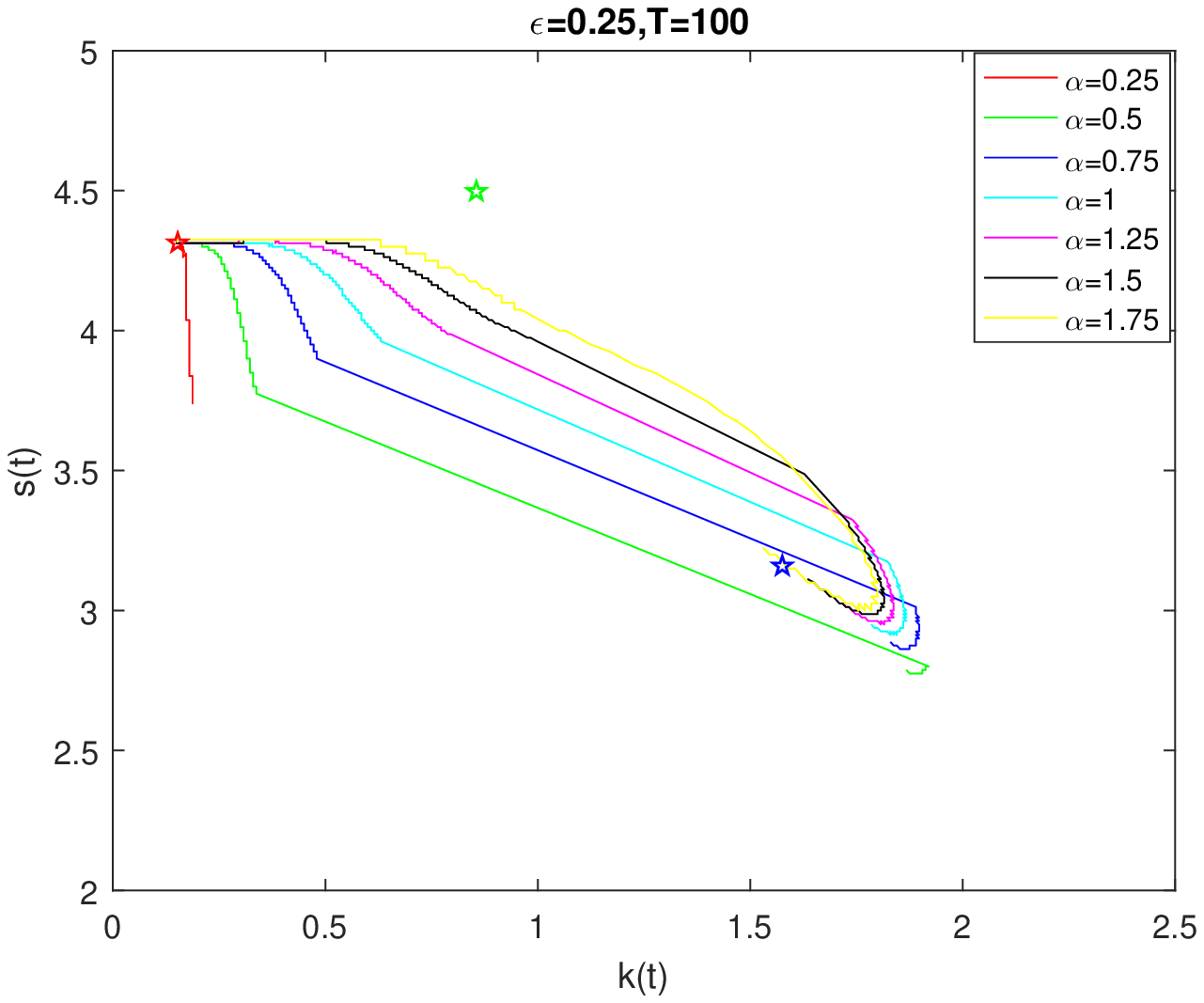}}
\end{minipage}
\hfill
\begin{minipage}[]{0.5 \textwidth}
 \leftline{~~~~~~~\tiny\textbf{(c2)}}
\centerline{\includegraphics[width=6cm,height=4cm]{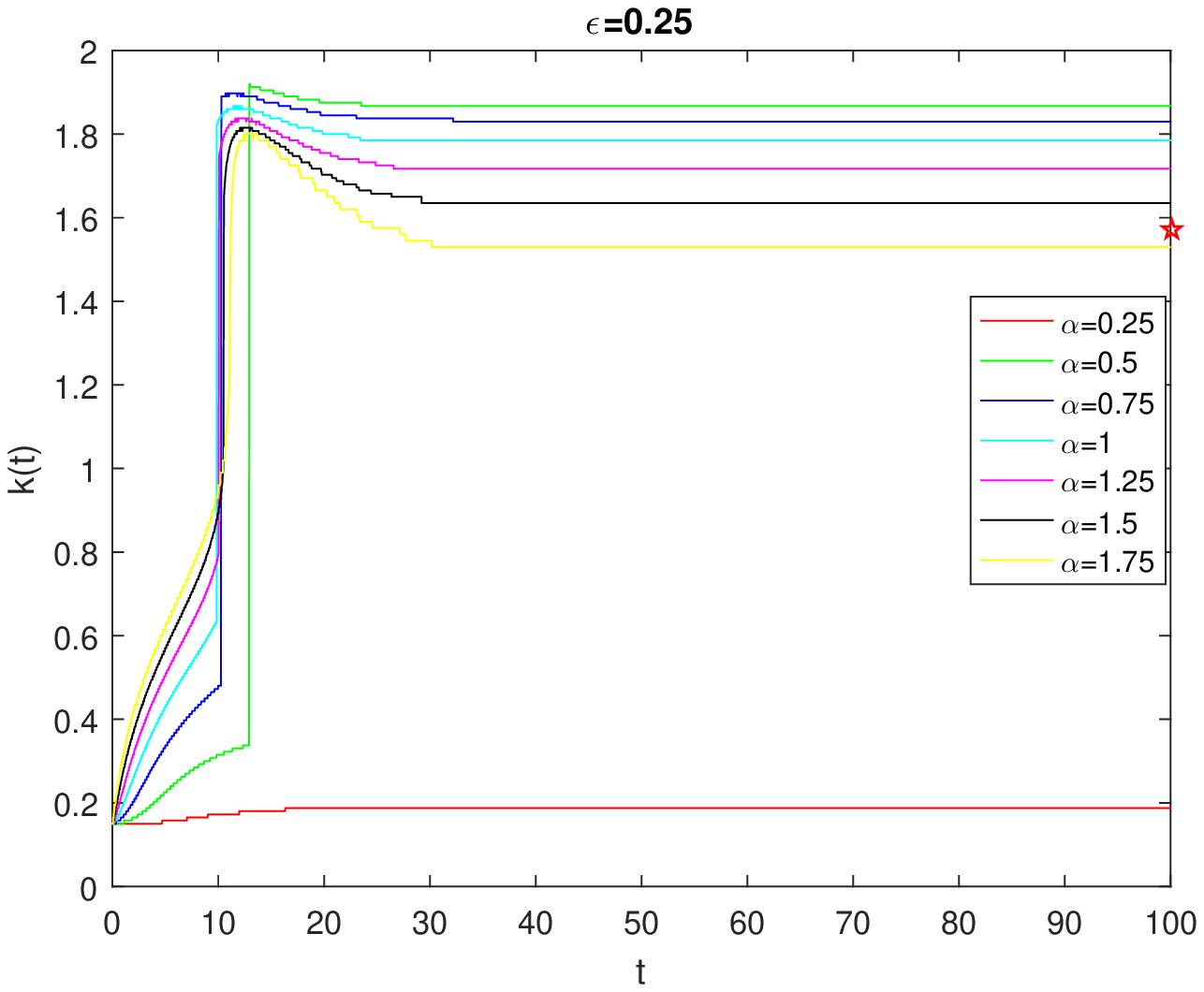}}
\end{minipage}

\begin{minipage}[]{0.5 \textwidth}
 \leftline{~~~~~~~\tiny\textbf{(d1)}}
\centerline{\includegraphics[width=6cm,height=4cm]{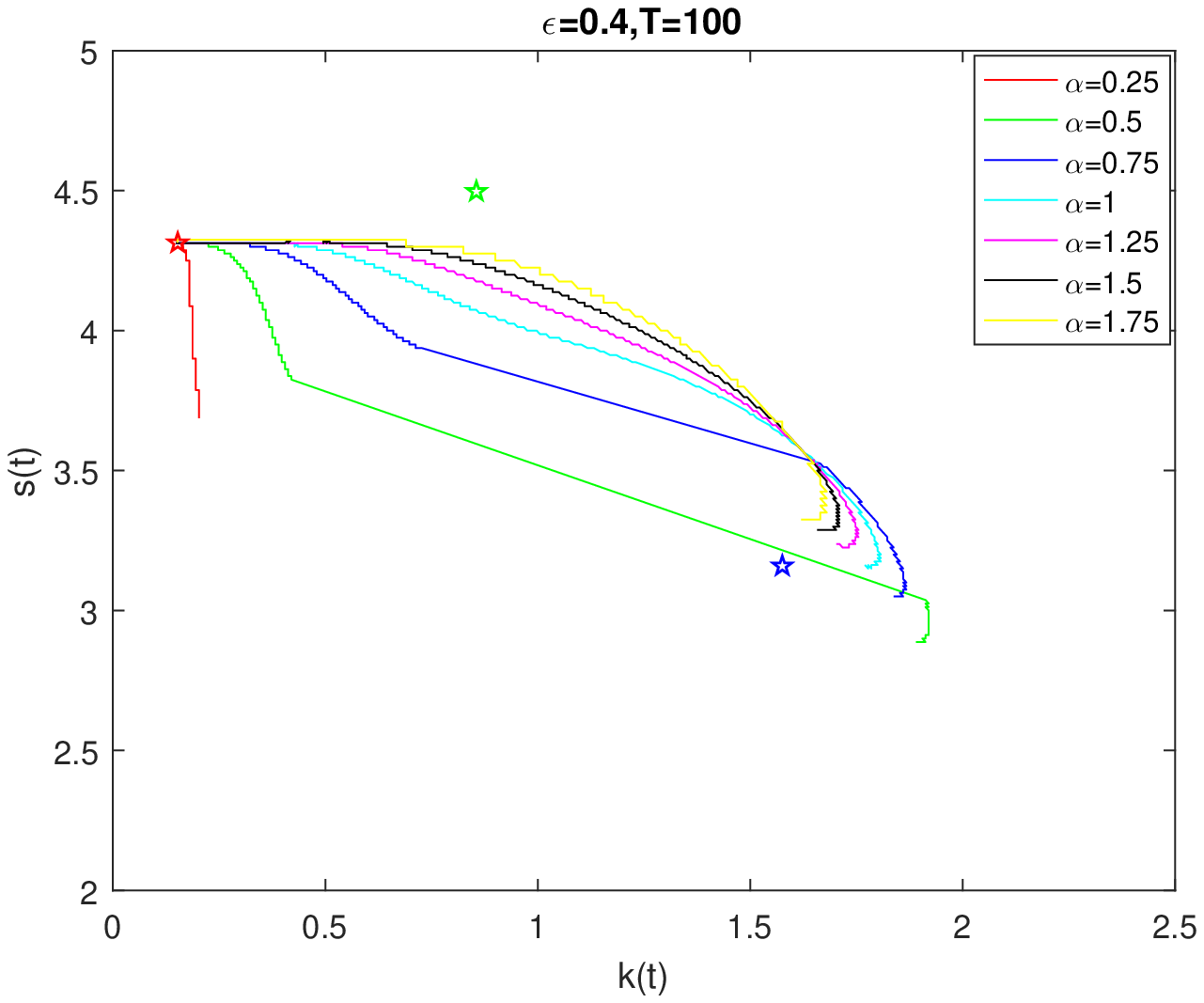}}
\end{minipage}
\hfill
\begin{minipage}[]{0.5 \textwidth}
 \leftline{~~~~~~~\tiny\textbf{(d2)}}
\centerline{\includegraphics[width=6cm,height=4cm]{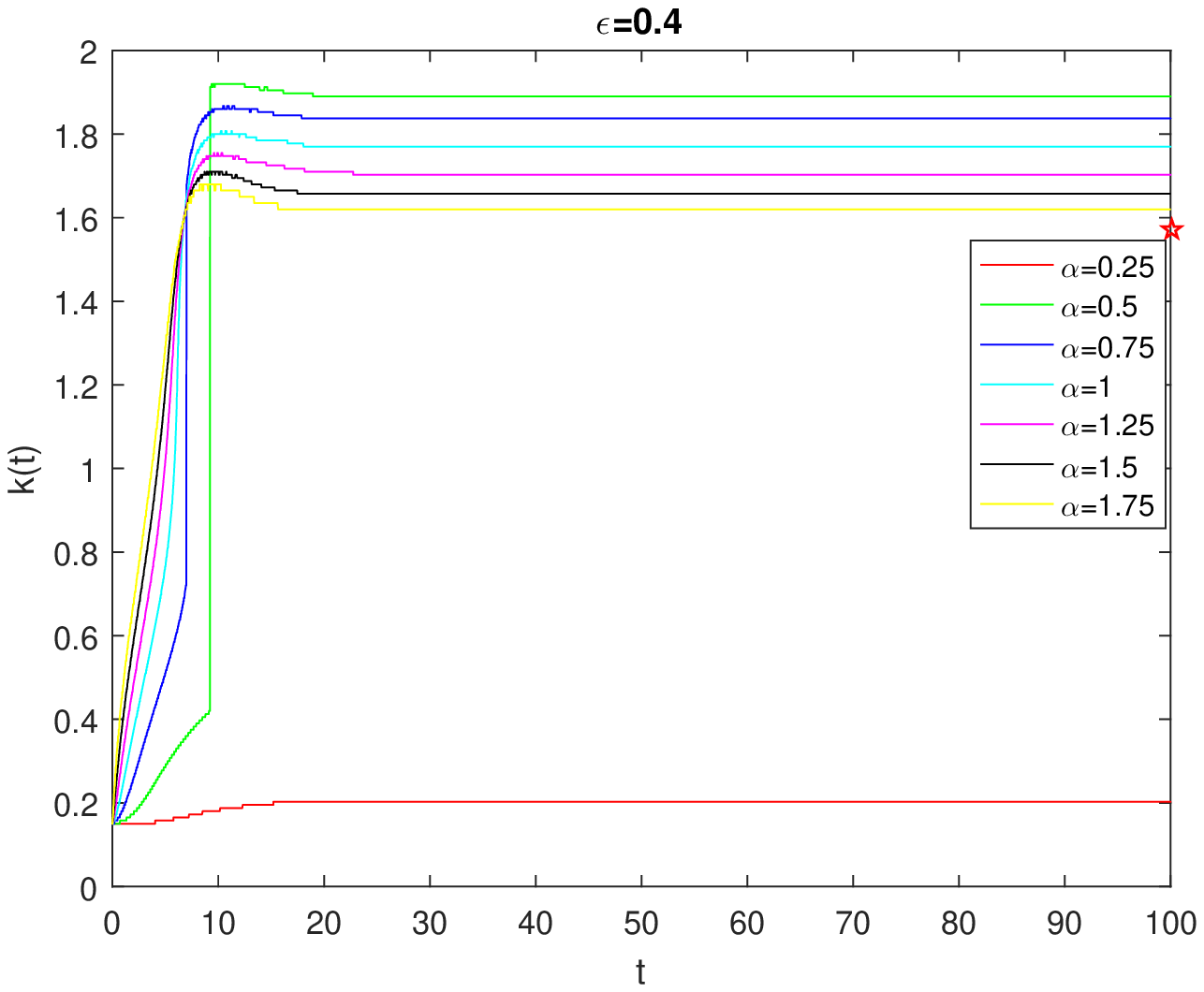}}
\end{minipage}

\vfill
\caption{\textbf{  Effect of $\alpha$ and $\epsilon$ for the most probable trajectories (left panel) and the corresponding time series of the ComK protein (right panel). } The most probable trajectories: (a1) $\epsilon$=0.1, (b1) $\epsilon$=0.15, (c1) $\epsilon$=0.25,  (d1) $\epsilon$=0.4. Time series of the ComK protein: (a2) $\epsilon$=0.1, (b2) $\epsilon$=0.15, (c2) $\epsilon$=0.25, (d2) $\epsilon$=0.4.}
\label{f4}
\end{figure}
\newpage
\textbf{Most probable trajectories and corresponding time series of the ComK protein for various noise intensity $\epsilon$ and non-Gaussianity index $\alpha$.} The left plots of FIG. \ref{f4} show the impacts of $\epsilon$ and $\alpha$ on the most probable trajectories. Note that the red star is the low concentration state, the green star denotes the saddle point and the blue star is the high concentration state. The right plots of FIG. \ref{f4} show the impacts of $\epsilon$ and $\alpha$ on the corresponding time series of ComK and the red star is the high concentration state of ComK protein.

When $\epsilon$ is small, the most probable trajectories stay around the low concentration state (FIGs. \ref{f4} (a1) and (a2) ). Hence, a small noise intensity is not in favor of the expression of \emph{comK} gene. When $\epsilon=0.15$ (see FIGs. \ref{f4} (b1) and (b2) ), the most probable trajectories first get near the unstable point, and then jump to the high concentration region. As time goes on, the most probable trajectories transit back to the low concentration region which means that the high concentration region is less stable or robust than the low concentration region in this case ($\epsilon=0.15$). With the increase of $\epsilon$, trajectories with $\alpha>0.25$ escape to the high ComK protein concentration region. But trajectories with $\alpha=0.25$ can not escape to the high concentration region by $T=100$. Note that when $\alpha=0.25$, the L\'evy noise has larger jumps but lower jump frequencies. FIGs. \ref{f4} (c1) and (d1) exhibit that when $\epsilon=0.25$ and $0.4$, trajectories  with $0.25<\alpha<2$ lead the ComK protein  evolve to the high concentration region, and in this case, the jump frequencies are higher than {\color{red}for} L\'evy noise with $\alpha=0.25$.

Furthermore, we see that when the most probable trajectories get close to the saddle point  \\
$(0.8568,4.4938)$, they will reach the high concentration soon. This is in accordance with the observation that the saddle point repulses the most probable trajectories, while the high concentration stable state appears to attract them. If the most probable trajectories pass the saddle point, then they arrive at the high concentration. We define the time spending from the low concentration state to the saddle point as the tipping time.

When the noise intensity is large, the tipping time becomes short (see FIGs. \ref{f4} (a2), (b2), (c2) and (d2) ). In the following, we will discuss the tipping time for different $\epsilon$ and $\alpha$.
\begin{figure}[H]
\begin{minipage}{0.48\linewidth}
\leftline{~~~~~~~\tiny\textbf{(a)}}
\centerline{\includegraphics[height = 6cm, width = 7cm]{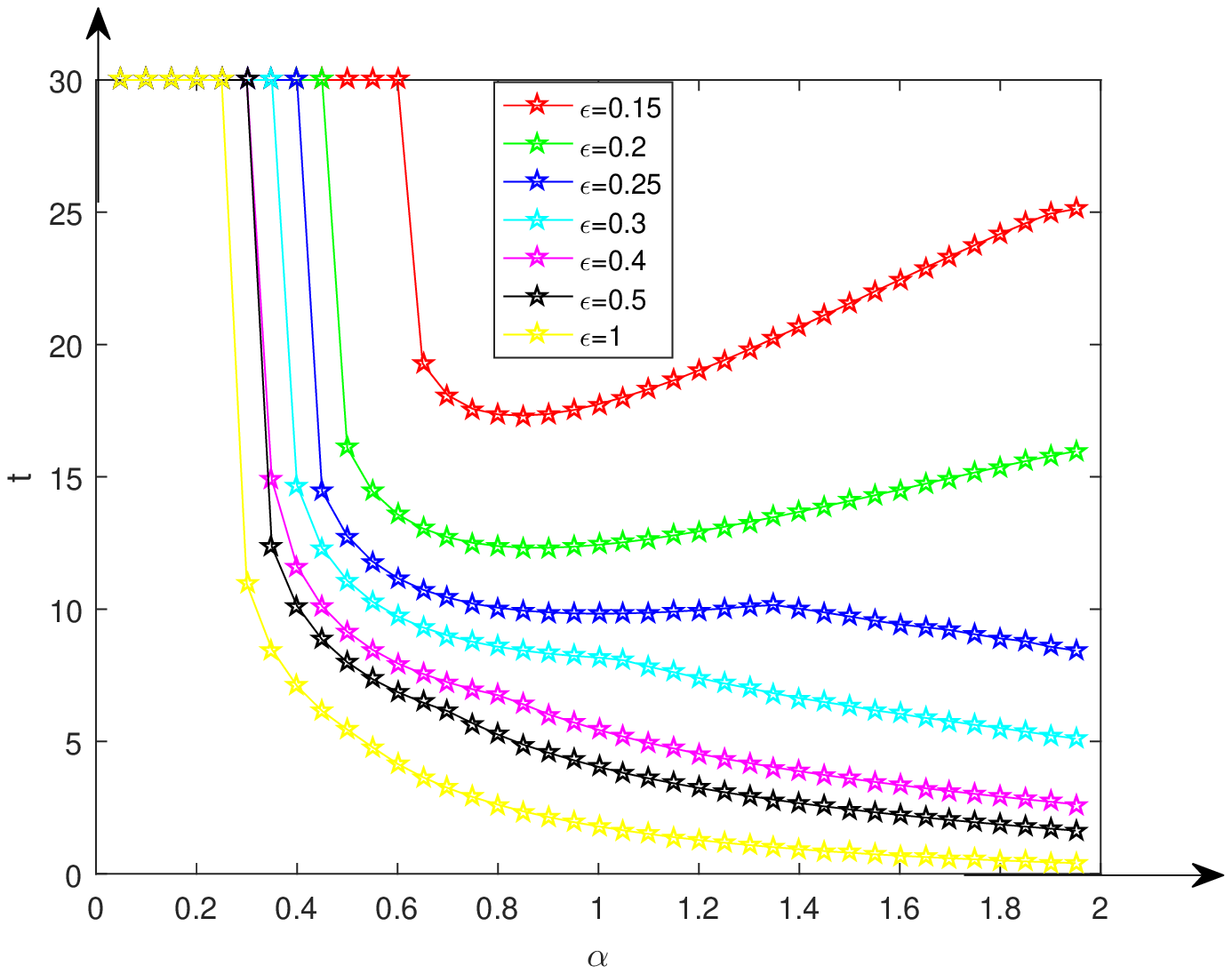}}
\end{minipage}
\hfill
\begin{minipage}{0.48\linewidth}
\leftline{~~~~~~~\tiny\textbf{(b)}}
\centerline{\includegraphics[height = 6cm, width = 7cm]{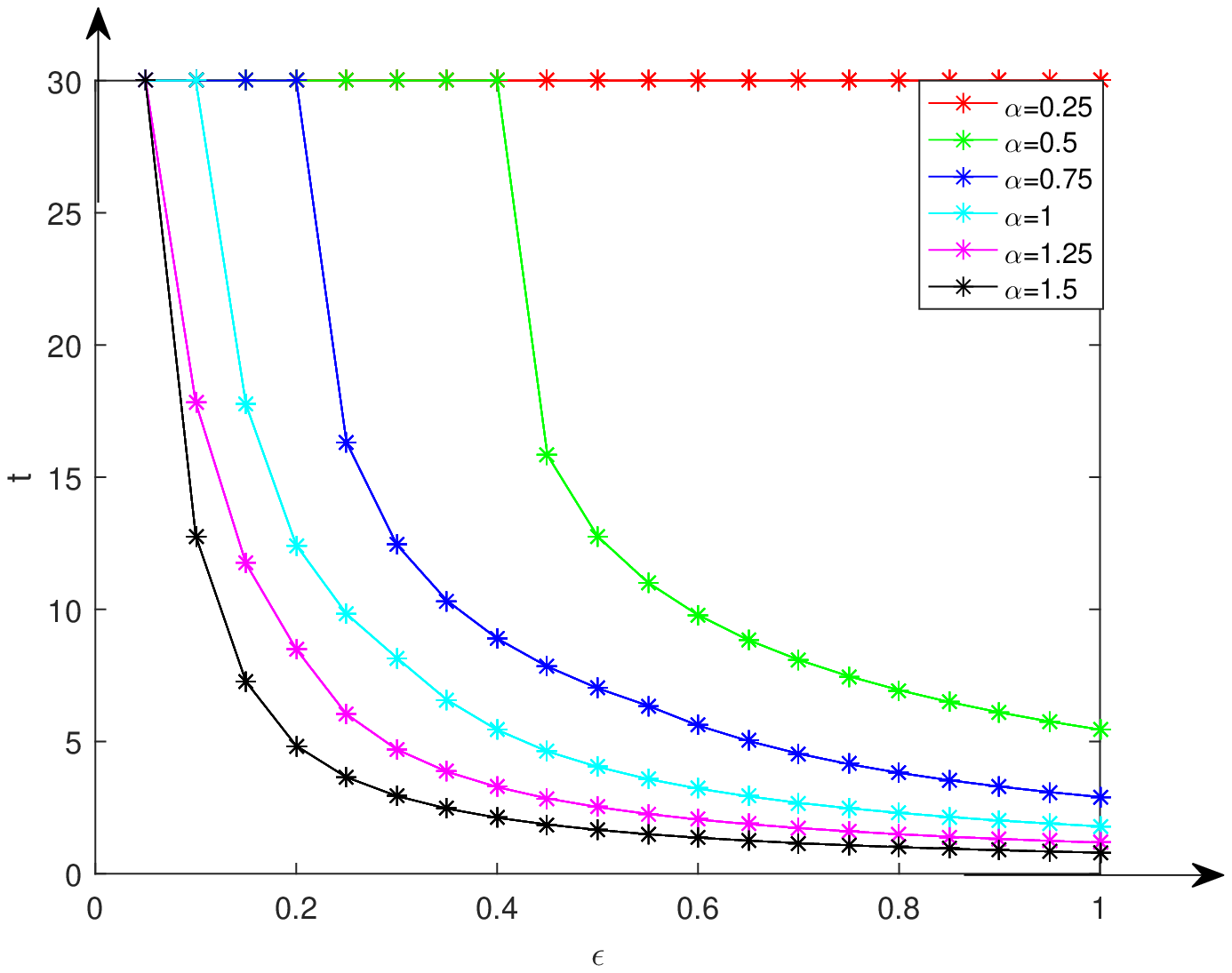}}
\end{minipage}
\caption{ \textbf{  Effect of $\alpha$ and $\epsilon$ for the tipping time from the low concentration to the high concentration.}   (a) Tipping time from the low concentration to the high concentration as a function of $\alpha$ with different $\epsilon$. (b) Same as (a) but as a function of $\epsilon$ with different $\alpha$.  }
  \label{f7}
\end{figure}

\textbf{Tipping time with different noise intensity $\epsilon$ and non-Gaussianity index $\alpha$.}
FIG. \ref{f7} shows the tipping time from the low concentration region to high concentration region with different $\epsilon$ and $\alpha$. It is worth pointing out that the most probable trajectories can not reach the high concentration region when the tipping time is up to 30.

As shown in FIG. \ref{f7} (a), when $\epsilon<0.2$, the tipping time first decreases and then increases as $\alpha$ increases. This means that $\alpha\approx 0.8$ with larger jumps is beneficial to making the expression of the \emph{comK} gene. When $\epsilon>0.2$, the tipping time decreases as $\alpha$ increases, i.e. large $\alpha$ (smaller jumps with higher frequency) helps the expression of the \emph{comK} gene. In this case, we may  infer that the transition will occur at earlier times for larger $\alpha$. These results coincide with results in \cite{huicfxuy}.

FIG. \ref{f7} (b) presents the tipping time as a function of $\epsilon$. The transition does not occur for $\alpha=0.25,$ due to its low jump frequency. For $\alpha>0.25,$ the tipping time decreases as the noise intensity increases. Except for $\alpha=0.25$, the tipping time decreases as $\alpha$ increases at a given noise intensity, because the jump frequency increases. In summary, a larger noise intensity with a larger non-Gaussianity index will play a more positive role in the expression of the \emph{comK} gene and a relative small noise intensity with non-Gaussianity index closing to $1$ will play a positive role in the expression of the \emph{comK} gene.
\begin{figure}[H]
\centering
\includegraphics[width=0.5\linewidth]{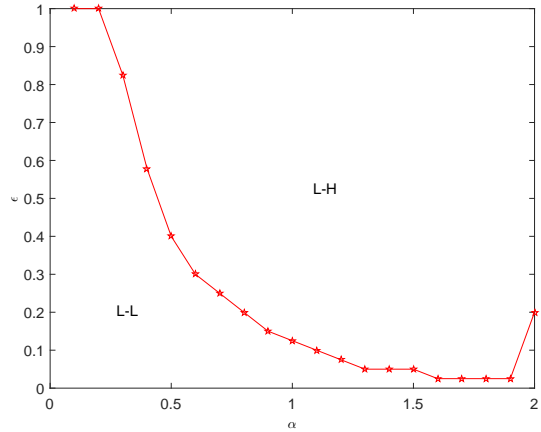}
\caption{
\textbf{Choice of combination of $\epsilon$ and $\alpha$ such that the ComK protein evolve from the low concentration to the high concentration.} L-L: the choice of combination of $\epsilon$ and $\alpha$ can not make the ComK protein evolve from the low concentration to the high concentration. L-H: the choice of combination of $\epsilon$ and $\alpha$ can make the ComK protein evolve from the low concentration to the high concentration.   }
\label{f5}
\end{figure}

\textbf{Combination of the noise intensity $\epsilon$ and the non-Gaussianity index $\alpha$ control the competence development.} Now we discuss which combination of $\alpha$ and $\epsilon$ can make the ComK protein {\color{red}efficiently} transit to the high concentration, i.e. the competence develops. It is shown in FIG. \ref{f5}, that the region $L-L$ below the red curve can not make a ComK protein transit from the low concentration to the high one, but that the region $L-H$ above the red curve can induce the ComK protein transit to the high concentration. This figure provides clearly how to choose $\epsilon$ and $\alpha$ to develop the competence. Larger $\alpha$ with high frequency jumps is propitious to the expression of the \emph{comK} gene. Small $\epsilon$ is not beneficial to the expression of the \emph{comK} gene.

\begin{figure}[H]
\begin{subfigure}[b]{0.4\textwidth}
\leftline{~~~~~~~\tiny\textbf{(a1)}}
        \includegraphics[width=\textwidth]{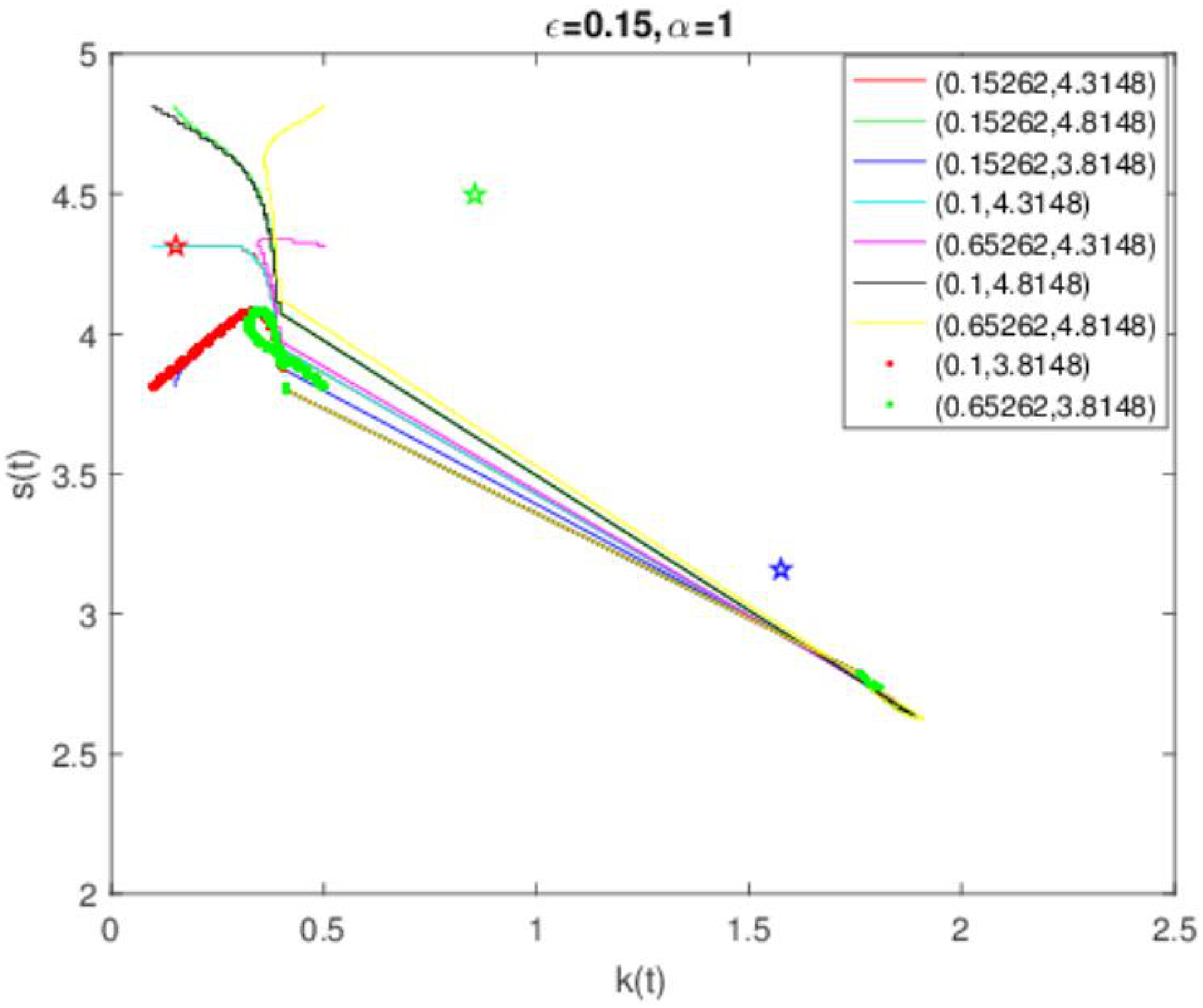}
    \end{subfigure}
    \hspace{3.0cm}
   \begin{subfigure}[b]{0.4\textwidth}
   \leftline{~~~~~~~\tiny\textbf{(b1)}}
        \includegraphics[width=\textwidth]{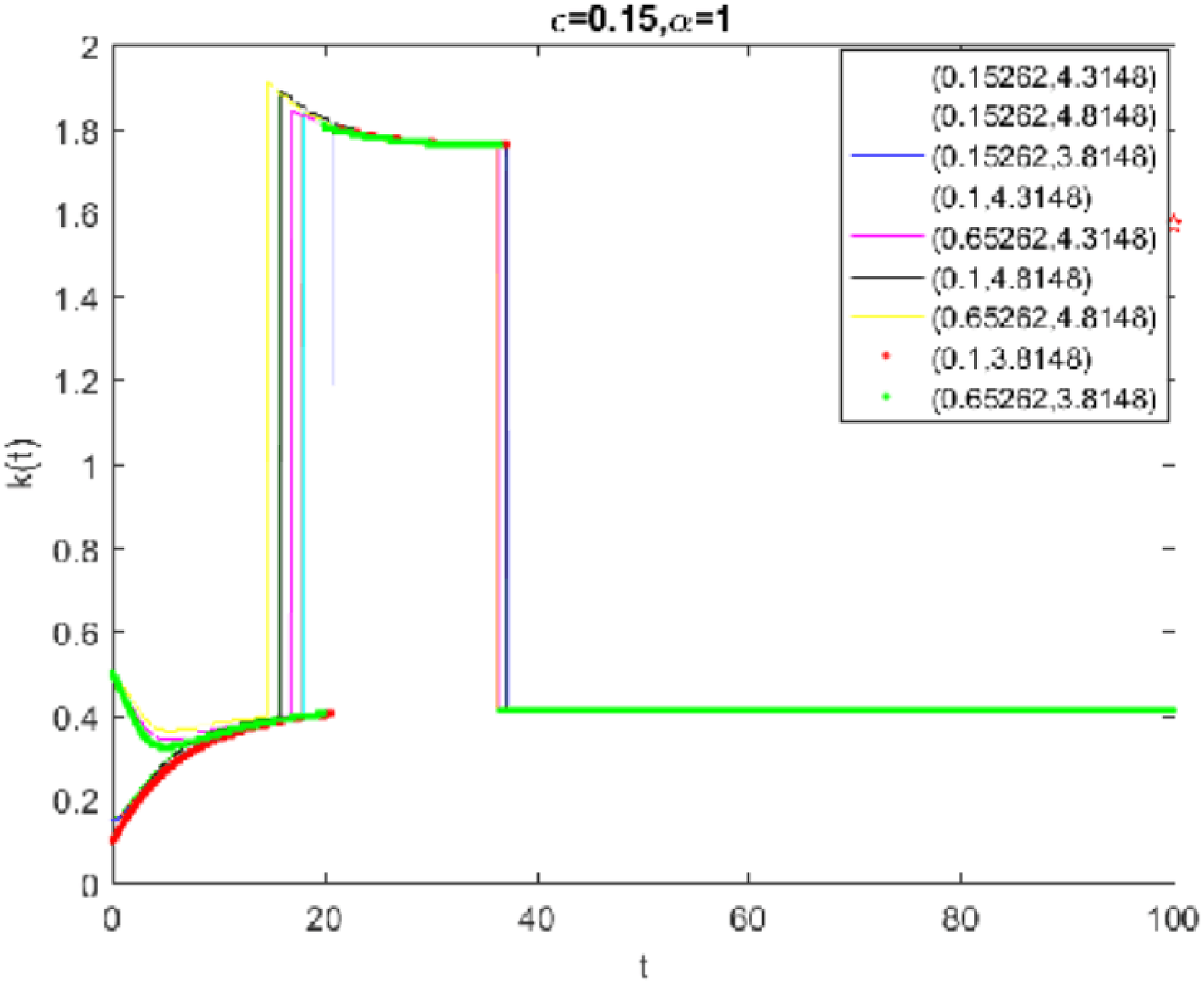}
    \end{subfigure}
    \begin{subfigure}[b]{0.4\textwidth}
    \leftline{~~~~~~~\tiny\textbf{(a2)}}
        \includegraphics[width=\textwidth]{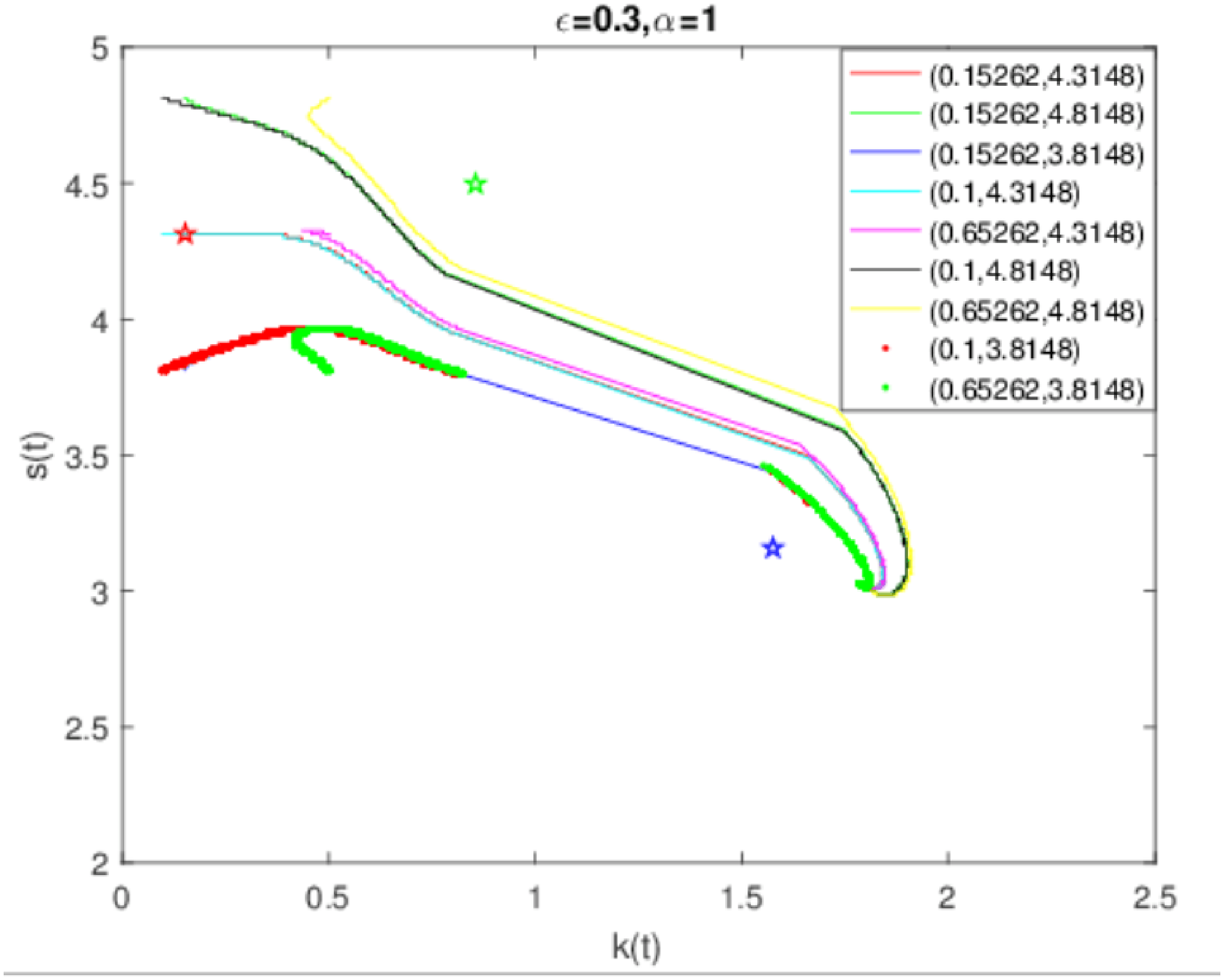}
    \end{subfigure}
    \hspace{3.0cm}
   \begin{subfigure}[b]{0.4\textwidth}
   \leftline{~~~~~~~\tiny\textbf{(b2)}}
        \includegraphics[width=\textwidth]{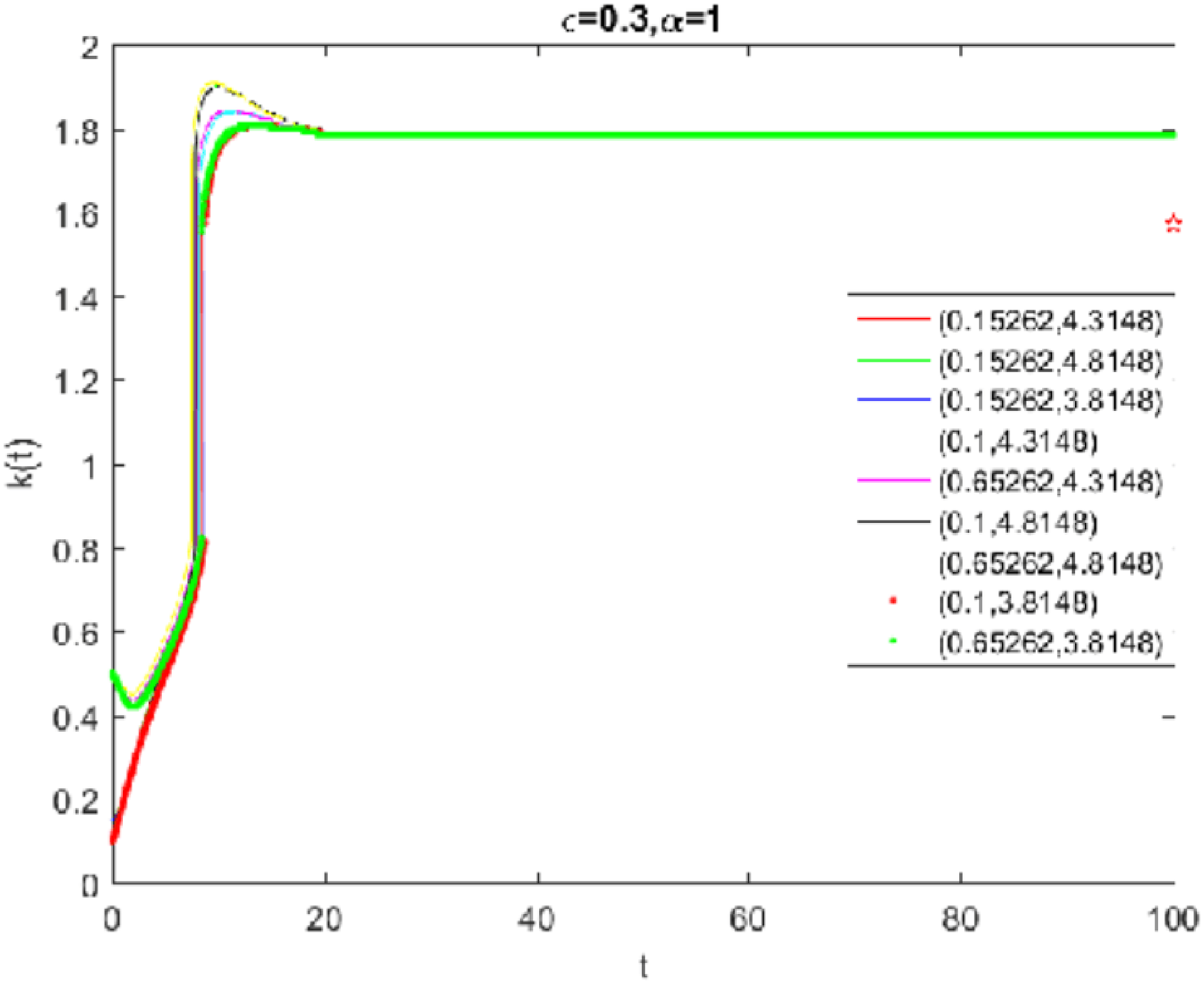}
    \end{subfigure}
\caption{ \textbf{  Effect of various initial concentrations for the most probable trajectories (left panel) and the corresponding time series of the ComK protein (right panel).}
(a) Most probable trajectories for different initial concentrations when $\epsilon$=0.15 and $\alpha$=1. (b) Time series of the ComK protein for different initial concentrations when $\epsilon$=0.15 and $\alpha$=1. (c) Most probable trajectories for different initial concentrations when $\epsilon$=0.3 and $\alpha$=1. (d) Time series of the ComK protein for different initial concentrations when $\epsilon$=0.3 and $\alpha$=1. }
\label{f8}
\end{figure}

   \textbf{Various initial concentrations evolve to a metastable state.} Here we study the effect of different initial concentrations near the nodal sink (the low ComK vegetative state) on the most probable trajectories (FIG. \ref{f8}). As in FIG. \ref{f4}, the red star denotes the nodal sink, the green star represents the saddle and the blue star is the spiral sink. As shown in FIGs. \ref{f8} (a1), (a2), when $\epsilon$ and $\alpha$ are fixed, we choose nine different initial concentrations surrounding the nodal sink, and find that the most probable trajectories show a slight tendency to the nodal sink, then ultimately gather at a metastable state, i.e. reaching the equilibrium in the sense of stochasticity.  FIGs. \ref{f8} (b1), (b2) display the most probable time series of the ComK protein, from which we can see that the ComK protein tends to the metastable state as time goes on. Interestingly, as shown in FIG. \ref{f8} (b1), when $\epsilon$ is relatively small ($\epsilon=0.15$), the most probable trajectories of the nine different initial concentrations first attain the competence region (the high ComK concentration), then they return to the vegetative region (the low ComK concentration), and finally gather at a low metastable state. However, when $\epsilon$ becomes larger ($\epsilon=0.3$), the most probable trajectories of the nine different initial concentrations reach the competence region (the high ComK concentration), and gather at a relative high metastable state finally. The differences among the nine most probable trajectories become narrow initially, then spread relatively strong, and eventually they converge at the high metastable state.

\begin{figure}[H]
\centering
\includegraphics[width=0.5\linewidth]{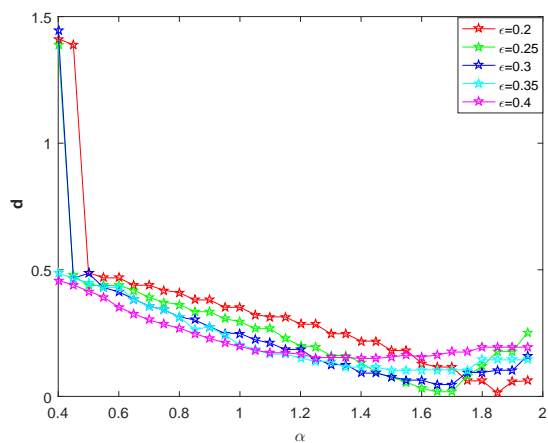}
\caption{Distance between the deterministic competence state and the metastable state for different $\alpha$ and $\epsilon$.}
\label{f9}
\end{figure}

   \textbf{The distance between the deterministic competence state and metastable state with different $\alpha$ and $\epsilon$.} As shown in FIGs. \ref{f8} (a1) and (a2), the most probable trajectories of different initial concentrations gather at different metastable states under the influence of L\'evy noise with different $\alpha$ and $\epsilon$. That is, the distance between the deterministic competence state and the metastable state varies with different $\alpha$ and $\epsilon$, from which we know that different L\'evy noise induce various effects on the behavior of dynamical systems. Now we discuss which $\alpha$ and $\epsilon$ can make the distance smallest, i.e. acting most efficiently. We consider the distance between the deterministic competence state $(k_{+},s_{+})$ and the metastable state $(k(T),s(T))$ as $\textbf{d}=\sqrt{(k_{+}-k(T))^2+(s_+-s(T))^2}$. It is shown in FIG. \ref{f9}, that $\textbf{d}$ decreases firstly and then increases with the increase of $\alpha$ when the noise intensity is fixed. When the noise intensity becomes stronger, the minimum of $\textbf{d}$ is bigger. When $\alpha=1.85$ and $\epsilon=0.2$, $\textbf{d}$ is the smallest. Hence, L\'{e}vy noise with $\alpha=1.85$ and $\epsilon=0.2$ have the smallest effect on the deterministic competence state.

\section{Conclusion}\label{Conclusion}
We have studied the most probable trajectories of a two-dimensional system (MeKS network) under the influence of non-Gaussian stable L\'{e}vy noise. On the basis of the nonlocal Fokker-Planck equation, we have described the most probable trajectories of the MeKS network from the low ComK protein concentration (vegetative state) to the high ComK protein concentration (competence state) under stable L\'{e}vy noise. It has been recently demonstrated that the non-Gaussian stable L\'{e}vy noise is more appropriate to describe the burst-like phenomenon in synthesis process of biochemical systems. Therefore, we here study the most probable trajectories of the MeKS network under the influence of stable L\'{e}vy noise.

The most probable trajectories of the MeKS network influenced by L\'evy noise with various noise intensity $\epsilon$ and non-Gaussianity index $\alpha$ have been investigated. We have found that for fixed $\alpha>0.25,$ the ComK protein stays at the low concentration region for small noise intensity, and it begins to transit to the high concentration region with the increase of $\epsilon$. Also, the tipping time from the low concentration region to the high concentration region as a function of $\alpha$ with various $\epsilon$ has been analysed. We have discovered that for fixed $\epsilon$, the tipping time firstly decreases and then increases as $\alpha$ increases. Thus, we can determine an optimal combination of $\alpha$ and $\epsilon$ to make the tipping time shortest. In addition, combinations of parameters of $\alpha$ and $\epsilon$ which generate a ComK protein transit from the low concentration to the high concentration have also been revealed. Furthermore, we have selected different initial concentrations around the low ComK protein concentration to examine the most probable trajectories of the MeKS network, to discover that the most probable trajectories gather at a metastable state. L\'evy motion with a small $\alpha$ ($\alpha\leq 0.25$) which has a lower jump frequency can not make the ComK protein transit to the competence region. For small $\epsilon$, L\'evy motion with a relative small $\alpha$, which has larger jumps, is more convenient to make the ComK protein transit to the competence region. However, for larger $\epsilon$, L\'evy motion with a larger $\alpha$ which has more small jumps is more convenient to transit to the competence region. Additionally, when $\alpha=1.85$ and $\epsilon=0.2$, the distance between the deterministic competence state and metastable state is smallest.

In short, we have utilised this geometric tool, the most probable trajectories, to visualise the dynamics of the MeKS networks, which offer an optimal path for the ComK protein escapes from the low concentration (vegetative region) to the high
concentration (competence region).

\section*{Appendixes}\label{Appendixes}
\subsection*{\textbf{Symmetric $\alpha$-stable L\'{e}vy motions}}\label{A}
   A stochastic process $\{ L_t: t\geq 0\}$ defined on $\mathbb{R}^d$ is a symmetric $\alpha$-stable L\'{e}vy process with $0<\alpha <2$ if the following conditions are satisfied:\begin{enumerate}
\item[(1)] $L_0=0$, a.s.;
\item[(2)] For any choice of $n\geq 1$ and $0< t_0 < t_1<\cdot\cdot\cdot<t_n$, random variables $L_{t_0},~L_{t_1}-L_{t_0},~L_{t_2}-L_{t_1},\cdot\cdot\cdot,~L_{t_n}-L_{t_{n\!-\!1}}$ are independent (independent increments property);
\item[(3)] The distribution of $L_{s+t}-L_s$ does not depend on $s$ (temporal homogeneity or stationary increments property);
\item[(4)] $L_t$ has stochastically continuous sample paths, i.e., for every $s>0,~L_t \rightarrow L_s$ in probability, as $t\rightarrow s$.
\end{enumerate}

\subsection*{\textbf{Fokker-Planck equation}}\label{B}
  The generator $\mathscr{A}$ of the process $(k_t,s_t)$ is defined as (see \cite{Applebaum,Cai,wu17}),
\begin{align}
\mathscr{A}u(k,s)&= f_1u_k + f_2u_s+ \int_{\mathbb{R}^2 \setminus \{(0,0)\} } [u(k\!\!+\!\!k',  s\!\!+\!\!s')   \nonumber   \\
&~~~~~-u(k, s)] [\epsilon_k^\alpha \nu_\alpha(dk') \delta_0(ds') + \epsilon_s^\alpha \nu_\alpha(ds') \delta_0(dk')]  \nonumber\\
 &=f_1 u_k+\epsilon_k^\alpha\int_{\mathbb{R} \setminus\{0\}} [u(k\!\!+\!\!k',  s)-u(k, s)
     ] \nu_\alpha(dk')                 \nonumber  \\
     &~~~~~+ f_2 u_s +\epsilon_s^\alpha\int_{\mathbb{R} \setminus\{0\}} [u(k,  s\!\!+\!\!s')-u(k, s)
     ]  \nu_\alpha(ds'),
\end{align}
 where $\delta_0$ is the delta measure concentrated at $0$ with property
$ \int_{-\infty}^{\infty} g(z) \delta_0(dz) = g(0)$.

The Fokker-Planck equation for the distribution of the conditional probability
density \\$p(k,s,t) = p(k,s,  t|k_0, s_0, 0)$, i.e., the probability of the
process $(k_t,s_t)$ has value $(k,s)$ at time $t$ given it had value $(k_0,s_0)$ at
time $0$, is  $p_t=\mathscr{A}^* p$  given by \cite{oksendal, DuanBook}, where $\mathscr{A}^*$ is
the adjoint operator for $\mathscr{A}$.

The adjoint operator for $\mathscr{A}$ can be found out, as the integral part is a symmetric operator.
Thus we have the Fokker-Planck equation:
\begin{align} \label{fpe}
 p_t(k,s,t)
  &= - (f_1 p)_k   +
    \epsilon_{k}^{\alpha} \int_{\mathbb{R}^1 \setminus\{0\}} [p(k\!\!+\!\!k',  s, t)-p(k, s, t)] \nu_\alpha(dk')                 \nonumber  \\
     &~~~~- (f_2 p)_s +
     \epsilon_{s}^{\alpha}\int_{\mathbb{R}^1 \setminus\{0\}} [p(k,  s\!\!+\!\!s', t)-p(k, s, t)]  \nu_\alpha(ds'),
\end{align}
with the initial condition: $p(k,s,0)=\delta(k-k_0,s-s_0)$.
\subsection*{\textbf{Numerical method}}\label{C}
   In the following, we present the numerical algorithms the case of $D=(a,b)\times(c,d)$ with the
absorbing condition. Noting that the absorbing condition dictates that
$p$ vanishes outside $D$, we can simplify Eq.~(\ref{fpe}) as
\begin{align}
 p_t(k,s,t) &=- (f_1p)_k   - (f_2p)_s    \nonumber  \\
     &-
 \frac{\epsilon_{k}^{\alpha} C_\alpha}{\alpha}\left[\frac{1}{(k\!-\!a)^\alpha}+\frac{1}{(b\!-\!k)^\alpha}\right]p(k,s,t)
 +\epsilon_{k}^{\alpha}C_\alpha \int_{a\!-\!k}^{b\!-\!k} \frac{p(k\!+\!k',  s, t) \!-\! p(k,  s, t)}{|k'|^{1\!+\!\alpha}}\; {\rm d}k'
  \nonumber  \\
     &-
 \frac{\epsilon_{s}^{\alpha}C_\alpha}{\alpha} \left[\frac{1}{(s\!-\!c)^\alpha}+\frac{1}{(d\!-\!s)^\alpha}\right] p(k,s,t)
 \!+\! \epsilon_{s}^{\alpha} C_\alpha \int_{c\!-\!y}^{d\!-\!y} \frac{p(k,  s\!+\!s', t) \!-\! p(k,  s, t)}{|s'|^{1\!+\!\alpha}}\; {\rm d}s' ,
\label{fpe1Dn3}
\end{align}
for $(k,s) \in D $; and $p(k,s,t)=0$ for $(k,s) \notin  D$.

    Set $v=\frac{2}{b\!-\!a}(k\!-\!a)\!-\!1$, $w=\frac{2}{d\!-\!c}(s\!-\!c)\!-\!1$, $p(\frac{b\!-\!a}{2}v\!-\!\frac{a\!-\!b}{2},
\frac{d\!-\!c}{2}w\!-\!\frac{c\!-\!d}{2},t)=P(v,w,t)$, i.e., $k=\frac{b\!-\!a}{2}v\!-\!\frac{a\!-\!b}{2}$, $s=\frac{d\!-\!c}{2}w\!-\!\frac{c\!-\!d}{2}$, then we get
\begin{align}
P_t(v,w,t) &=-\frac{2}{b-a}(f_1 P)_v   -\frac{2}{d-c} (f_2 P)_w \nonumber  \\
&
 - \frac{ C_\alpha}{\alpha} (\frac{2\epsilon_{k}}{b\!-\!a})^{\alpha} \left[\frac{1}{(1\!+\!v)^\alpha}\!+\!\frac{1}{(1\!-\!v)^\alpha}\right] P\!+\! C_\alpha (\frac{2\epsilon_{k}}{b\!-\!a})^{\alpha}\int_{-1\!-\!v}^{1\!-\!v} \frac{P(v\!+\!v',w,t) \!-\! P(v,w,t)}{|v'|^{1\!+\!\alpha}}\; {\rm d}v'\nonumber  \\
&
 \!-\!\frac{C_\alpha}{\alpha} (\frac{2\epsilon_{s}}{d\!\!-\!\!c})^{\alpha}\!\!\left[\frac{1}{(1\!\!+\!\!w)^\alpha}\!+\!\frac{1}{(1\!\!-\!\!w)^\alpha}\right] P\!+\! C_\alpha (\frac{2\epsilon_{s}}{d\!-\!c})^{\alpha}\!\!\int_{-1\!-\!w}^{1\!-\!w} \frac{P(v,w\!+\!w',t) \!-\! P(v,w,t)}{|w'|^{1\!+\!\alpha}}\; {\rm d}w'.
\label{fpe1Dn4}
\end{align}
for $(v,w) \in D'=(-1,1)\times(-1,1) $; and $P(v,w,t)=0$ for $(v,w) \notin  D'$.

Then, we use a numerical method to discretize the nonlocal Fokker-Planck equation \eqref{fpe1Dn4}. For the spatial direction, the advection term $- (f_1P)_x$ and $- (f_2P)_y$ are discretized by the third-order WENO method given in \cite{Jiang} and the singular integral term are used a modified trapezoidal rule to approximate\cite{Gaoting16,Gaoting14}. We divide the interval $[-2,2]\times[-2,2]$ in space into $(4I)^{2}$ sub-intervals and define $v_i=ih$, $w_j=jh$ for $-2I\leq i,j \leq 2I$, where $h=1/I$. Denoting the numerical solution of $P$ at $(v_i,w_j,t)$ by $P_{i,j}$, we obtain the semi-discrete equation
\begin{equation}
  \begin{split}
  \dfrac{{\rm d}P_{i,j}}{{\rm d}t} &=
   - \frac{2}{b-a}[(f_1P)_{v,i}^+ +  (f_1P)_{v,i}^-]-\frac{2}{d-c}[(f_2P)_{w,j}^+ +  (f_2P)_{w,j}^-] \nonumber  \\
   &-  \frac{ C_\alpha}{\alpha} (\frac{2\epsilon_{k}}{b\!-\!a})^{\alpha} \left[\frac{1}{(1+v_i)^\alpha}+\frac{1}{(1-v_i)^\alpha}\right] P_{i,j}
   -  \frac{ C_\alpha}{\alpha} (\frac{2\epsilon_{s}}{d\!-\!c})^{\alpha} \left[\frac{1}{(1+w_j)^\alpha}+\frac{1}{(1-w_j)^\alpha}\right] P_{i,j} \nonumber  \\
   & +C_{hx} \frac{P_{i-1,j} - 2P_{i,j} + P_{i+1,j}}{h^2}+
  C_{hy} \frac{P_{i,j-1} - 2P_{i,j} + P_{i,j+1}}{h^2}  \nonumber  \\
  & + C_\alpha (\frac{2\epsilon_{k}}{b\!-\!a})^{\alpha} h \sum^{I-i}_{k_1=-I-i,k_1\neq 0}\!\!\!\!\!\!\!\!\! \;
    {\frac{P_{i+k_1,j} - P_{i,j}}{|{v_k}_1|^{1+\alpha}} }
    + C_\alpha (\frac{2\epsilon_{s}}{d\!-\!c})^{\alpha} h \sum^{I-j}_{k_2=-I-j,k_2\neq 0}\!\!\!\!\!\!\!\!\! \;
    {\frac{P_{i,j+k_2} - P_{i,j}}{|{w_k}_2|^{1+\alpha}} },
  \end{split}
 \label{nm1D3}
\end{equation}
for $0<\alpha<2$ and $i,j = -I+1, \cdots, -2,-1,0,1,2, \cdots, I-1$, where
the constant $C_{hx}\!=\! -  C_\alpha (\frac{2\epsilon_{k}}{b\!-\!a})^{\alpha}
\zeta(\alpha-1) h^{2-\alpha}$ and $C_{hy}\!=\! - C_\alpha (\frac{2\epsilon_{s}}{d\!-\!c})^{\alpha}
\zeta(\alpha-1) h^{2-\alpha}$, $\zeta$ is the Riemann zeta function, and the $\pm$ superscripts denote the global Lax-Friedrichs flux splitting defined as $(f_kP)^{\pm}_{i,j} =
\frac{1}{2}(f_kP_{i,j} \pm \alpha_k P_{i,j})$ with $\alpha_k=\max |f_k(x,y)|$,~$k=1,2$. The
absorbing condition requires that the values of $P_{i,j}=0$ if the index $|i|\geq I$ or $|j|\geq I$.

For time discretization, we use a third-order total variation
diminishing Runge-Kutta method provided in \cite{Shu}. In
particular, for the ordinary differential equation $\dfrac{{\rm
d}P}{{\rm d}t}= R(P)$, the method can be written as
\begin{equation}
\begin{split}
P^{(1)} &= P^n + \Delta t R(P^n),\\
P^{(2)} &= \frac{3}{4}P^n + \frac{1}{4} P^{(1)} +\frac{1}{4}\Delta t
R(P^{(1)}),\\
P^{n+1} &= \frac{1}{3} P^n + \frac{2}{3} P^{(2)} + \frac{2}{3}
\Delta t R(P^{(2)}),
\end{split}
\label{eq.rk3}
\end{equation}
where $P^n$ denotes the numerical solution of $P$ at time $t=t_n$.

\section*{Acknowledgements}
We would like to thank Dongfang Li, Ke Yin, Yayun Zheng, Xiujun Cheng and Rui Cai for helpful discussions.

\end{document}